# Simulation of Simplicity:

# A Technique to Cope with Degenerate Cases in Geometric Algorithms[1]

Herbert Edelsbrunner[2] and Ernst Peter Mücke[2]

**Abstract**

This paper describes a general-purpose programming technique, called the **S**imulation **o**f **S**implicity, which can be used to cope with degenerate input data for geometric algorithms. It relieves the programmer from the task to provide a consistent treatment for every single special case that can occur. The programs that use the technique tend to be considerably smaller and more robust than those that do not use it. We believe that this technique will become a standard tool in writing geometric software.

[1]Research of both authors was supported by Amoco Foundation Faculty Development Grant CS 1–6–44862. It was partially carried out while both authors were with the Institutes for Information Processing at the Technical University of Graz, Austria. The first author also acknowledges support by the National Science Foundation under grant CCR–8714565.

[2]Department of Computer Science, University of Illinois at Urbana-Champaign, 1304 West Springfield Avenue, Urbana, Illinois 61801, USA.

# 1 Introduction

This paper introduces a general technique that can be used to cope with degenerate cases encountered by computer programs. Consider, for example, a program that sorts an array of integers using a comparison as a primitive operation. A special, or degenerate, case occurs when the program attempts to decide which one of two equal numbers is smaller than the other. A typical way to resolve this tie is to pretend that the number with smaller index is smaller (assuming the integers are indexed, e.g., by their positions in an array). Or think of Kruskal's algorithm for constructing a minimum spanning tree of a weighted graph (see [AHU74]). At each step it chooses the shortest edge that can be added to the current collection of edges without creating a cycle. If this edge is not unique, then any one of the candidate edges is taken. The thus generated minimum spanning tree is therefore not unique unless we specify deterministic rules to break ties.

In both problems, sorting and constructing minimum spanning trees, the special cases are easily dealt with, partly because the ties can be broken arbitrarily without creating inconsistencies. The situation is usually far more complicated for geometric problems. Consider for example the following seemingly straightforward algorithm for the point-in-polygon problem which is sometimes called the *Parity Algorithm.*

- Let $r$ be the horizontal half-line whose left endpoint is the test point.
- Count the number of intersections between $r$ and the edges of the polygon. If that number is odd, then the test point lies within the polygon, and if the number is even, then it lies outside the polygon.

As pointed out in [Fo85], it is not a trivial matter to implement this algorithm, even if we assume that the test point does not lie on the boundary of the polygon. There are only two nondegenerate cases: Either the intersection between $r$ and an edge $e$ is empty or $r$ crosses $e$ (see Figure 1-I, (a) and (b)). There are, however, four degenerate cases (as illustrated in Figure 1-I, (c) through (f)) that have to be taken into account.

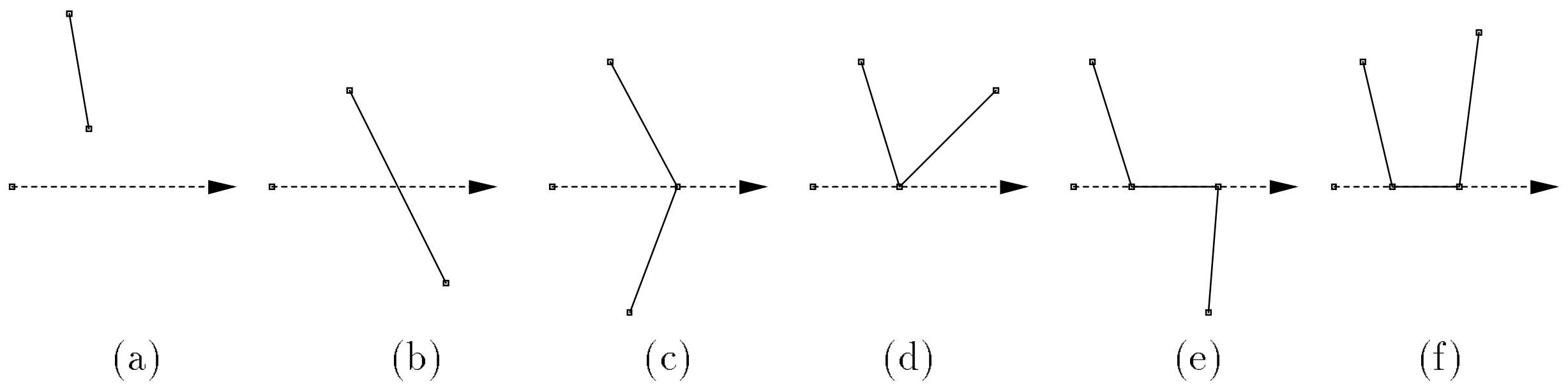

Figure 1-I: The different cases in the *Parity Algorithm.*

A correct answer is obtained if cases (c) and (e) are counted as one crossing and cases (d) and (f) are not counted at all. If we write the code for the above algorithm, we realize that a substantial amount of the effort is required to cover the four degenerated cases. Observe also that there are several seemingly plausible ways to treat the degenerate cases and that some of them lead to incorrect algorithms. We appeal to the imagination of the reader to envision the bizarre structure of degenerate cases one encounters in generalizing the point-in-polygon problem to three or higher dimensions. Another problem with a set of degenerate cases that is considerably richer than the one of the point-in-polygon problem is obtained if one intersects a polygon with a geometric object that is more complicated than a half-line.

When it comes to implementing geometric algorithms, degenerate cases are very costly, in particular if there are many such cases that have to be distinguished. This is caused by the positive correlation between the number of degenerate cases and a variety of factors that contribute to the overall cost of a piece of software. These factors include the length of the program, which, for itself, correlates positively with the amount of time required to write it, to debug it, and to maintain it. Of course, the degree of robustness of the program decreases with increasing complication. The correctness of a program relies on the consistent treatment of all different cases. In this context, it is worthwhile to mention that more efficient algorithms tend to be more complicated and also more sensible to slight inconsistencies in treating degenerate cases.

This paper presents a general technique, called **S**imulation **o**f **S**implicity (SoS), that can be used to cope with the problems mentioned above. Intuitively, it simulates a conceptual perturbation of the input data that eliminates all degeneracies. We hasten to mention that the perturbation is never ever computed — it is assumed to be arbitrarily small, although not vanishing, which is enough to simulate the nondegenerate topology. Another interpretation of the technique views it as a general way to break ties consistently. The tie-breaking part of the code appears in the lowest level of the algorithm, namely, in the procedures that implement the needed primitive operations. Different techniques following the same main approach have recently been suggested in [Ya87, Ya88]. A large part of this paper is devoted to demonstrating that the overhead in time caused by the use of the more elaborate primitive procedures required by SoS is negligible.

The outline of this paper is as follows: Section 2 presents the general idea of the technique and works out some guidelines needed to implement it effectively. Section 3 considers a class of problems for finite point sets that can be solved using a common set of geometric primitives. It also discusses how the perturbation influences the geometric primitives. Section 4 demonstrates efficient implementations of the primitive operations. In Section 5 we show that the geometric primitives introduced for point set problems can be used to solve a variety of other problems defined for polygons, hyperplanes, circles, spheres, and other geometric objects. Finally, in Section 6, we discuss the perturbation technique and its limitations.

## 2 SoS — the General Idea

Degeneracies occur with probability zero if we draw a finite number of geometric objects, each represented by a finite set of numbers from the (infinite) set of all such objects, provided there is no bound on the precision of the numbers used. In real-life computing this is not the case; that is, there is only a finite set of available numbers and thus a bound on the precision that can be achieved. As a consequence, we are doomed to work with degenerate data. On the other hand, even infinite precision does not guarantee the nonexistence of degeneracies. This section gives the general outline of a technique called the Simulation of Simplicity (SoS) — we use "simple" as a synonym for nondegenerate — which allows us to neglect degeneracies when we write programs. A similar but less elaborate method has been used to solve degenerate linear programs. This leads to the implementation of the simplex algorithm referred to as the "lexicographical method" (see [Ch52], [DOW55], [Da63], or [Ch83] for details). In computational geometry, this technique has been used in a couple of papers, including [Ed86] and [EW86], to avoid the otherwise necessary discussion of degenerate cases. This paper presents the theoretical foundations of SoS as well as details of its implementation.

The basic idea of SoS is to perturb the given objects slightly which amounts to changing the numbers that represent the objects; these numbers will be called the *coordinates* or the *parameters* of the objects. It is important that the perturbation is small enough so that it does not change the nondegenerate position of objects relative to each other. Coming up with such a perturbation is rather difficult and may require much higher precision than used for the original set of objects. For this reason, we perform the perturbation only symbolically by replacing each coordinate by a polynomial in $\varepsilon$. The polynomials will be chosen in such a way that the perturbed set goes towards the original set as $\varepsilon$ goes to zero. We will see that it is not important to know the exact value of $\varepsilon$ to perform the simulation; rather, it is sufficient to assume that $\varepsilon$ is positive and sufficiently small. Thus, it will be possible to use $\varepsilon$ as an indeterminant and to handle primitive operations symbolically.

The future user of SoS will neither have to be concerned with the role that $\varepsilon$ plays in the perturbation nor with the symbolic manipulation of polynomials. We may think of SoS as a package that provides the primitive operations needed for a certain computation. Ideally, the inside of these operations is hidden from the user who communicates with them like with an oracle. It turns out that a large number of geometric problems can be solved using a surprisingly small number of primitives. Some of these primitives will be discussed in the following three sections. This section continues to develop the general ideas on which SoS is based.

One of the goals of SoS is to perturb a set of objects such that all degeneracies disappear. A *degeneracy* is something that is not defined in general; its definition depends on the problem at hand. More specifically, it depends on the primitive operations used to solve the problem. For example, a primitive operation in the point-in-polygon algorithm described in the introduction tests the intersection of a horizontal half-line and a line segment. A degeneracy occurs if the half-line contains one or both endpoints of the line segment. A set of objects is now called *simple*, or *nondegenerate*, or *in general position*, if it does not contain any degeneracy. We thus define "simplicity" relative to the primitives used to solve a problem.

This paper considers only *topological primitives*, that is, operations that test some given input and classify it as one of a constant number of possible cases. This is in contrast to operations that compute new objects such as the intersection of a half-line and a line segment. In most programs, such an object serves only as an intermediate result anyway; but an intermediate result can as well be represented implicitly as a collection of pointers and a tag that tells us in what sense the objects identified by the pointers determine the (implicit) result. To simplify our discussion even further, we restrict our attention to primitives with three possible outcomes which we represent by $+1$, $0$, and $-1$, where $0$ indicates a degeneracy and $+1$ and $-1$ distinguish between the two nondegenerate cases. Tests that distinguish between more than two nondegenerate cases can be obtained by combining several ternary tests.

If we think of a primitive operation as a function $f$ that maps a high-dimensional point (whose coordinates describe the input objects) to $+1$, $0$, or $-1$, then $f^{-1}(0)$ represents the set of degenerate inputs. One requirement for this set is that its measure in this high-dimensional space is zero — otherwise, it is unreasonable to call its points degenerate. A set of $n$ objects, given by $d$ parameters each, can be thought of as a point in $nd$ dimensions. If $f$ takes $k < n$ objects as input, then $f^{-1}(0)$ is a surface of measure zero in $kd$-dimensional space. This surface defines another zero-measure surface in $nd$ dimensions which is obtained by embedding $f^{-1}(0)$ in the $kd$-dimensional subspace defined by the $k$ objects and extending it orthogonal to this subspace along the other coordinate axes. Other combinations of $k$ objects provide additional zero-measure surfaces that, altogether, decompose the $nd$-dimensional space into *faces* of various dimensions. A *cell* is an $nd$-dimensional face of this decomposition, and all points of a cell correspond to nondegenerate sets of objects. A degenerate set corresponds to a point $x$ in the union of the surfaces, denoted by $\mathcal{S}$. Since $\mathcal{S}$ has measure zero, every nonempty open ball around this point contains a point $y$ of some cell. Moving $x$ to $y$ corresponds now to perturbing the set of objects that $x$ corresponds to such that all degeneracies disappear. This shows that a perturbation to a nondegenerate set is always possible even if the amount of perturbation is severely limited. Recall that another requirement for the perturbation is that it does not change any nondegenerate subconfiguration. This means that we should not move $x$ across a surface it did not belong to initially. This can always be guaranteed if we choose the open ball small enough that it does not intersect any surface that does not contain the initial position of $x$.

To follow the forthcoming reasoning it is not necessary for the reader to understand the topology of the $nd$-dimensional space as indicated in the above paragraph. Nevertheless, this view of the problem sheds some light on the nature of degeneracy. It also explains why there is always a small enough perturbation that removes all degeneracies. Below, we discuss such perturbations more specifically and address a few questions concerning the efficient implementation of SoS.

Simplicity is simulated by applying a particular perturbation to a set $P = \{p_0, p_1, \ldots, p_{n-1}\}$ of $n$ geometric objects

$$p_i = (\pi_{i,1}, \pi_{i,2}, \ldots, \pi_{i,d}), \quad 0 \le i \le n-1,$$

each specified by $d$ parameters. It will be important that each object has a unique index between 0 and $n-1$. The objects are in arbitrary, and therefore not necessarily in simple, position. The perturbation of $P$ is realized by replacing each parameter by a polynomial in $\varepsilon$. We define

$$P(\varepsilon) = \{p_i(\varepsilon) = (\pi_{i,1}(\varepsilon), \pi_{i,2}(\varepsilon), \ldots, \pi_{i,d}(\varepsilon)) \,|\, 0 \le i \le n-1\},$$

where

$$\pi_{i,j}(\varepsilon) = \pi_{i,j} + \varepsilon(i,j) \quad \text{for} \quad 0 \le i \le n-1, \quad 1 \le j \le d,$$

and $\varepsilon(i,j)$ a polynomial in $\varepsilon$ that goes to zero when $\varepsilon$ goes to zero. We will refer to the new parameters $\pi_{i,j}(\varepsilon)$, the new objects $p_i(\varepsilon)$, and the new set $P(\varepsilon)$ as the *$\varepsilon$-expansions* of the original parameters $\pi_{i,j}$, the original objects $p_i$, and the original set $P$, respectively. The choice of the polynomials $\varepsilon(i,j)$ will be guided by three requirements SoS has to meet.

(a) $P(\varepsilon)$ must be simple if $\varepsilon > 0$ is sufficiently small.
(b) $P(\varepsilon)$ must retain all nondegenerate properties of the original set $P$.
(c) The computational overhead caused by simulating $P(\varepsilon)$ should be negligible.

As mentioned before, condition (b) is automatically met if $\varepsilon$ is small enough. To satisfy (a), it is sufficient to choose the $\varepsilon(i,j)$ such that there is no nonempty open interval $I$ with the property that $P(\varepsilon)$ is not simple if $\varepsilon \in I$. Think of $P$ as a point $x$ in $nd$ dimensions and let $x(\varepsilon)$ be the point that corresponds to $P(\varepsilon)$. The points $x(\varepsilon)$, $\varepsilon > 0$, form a one-dimensional curve $C$ in $nd$ dimensions. Thus, (a) is satisfied if $C \cap \mathcal{S}$ is a discrete set of points. (Recall that $\mathcal{S}$ represents all points in $nd$ dimensions that correspond to degenerate sets $P$.) In this topological setting, the phrase "$\varepsilon$ sufficiently small" gets a specific meaning: If $\varepsilon_0 > 0$ is the smallest value of $\varepsilon$ such that $x(\varepsilon_0) \in \mathcal{S}$, then $\varepsilon$ is sufficiently small if and only if $0 < \varepsilon < \varepsilon_0$. It is less clear how condition (c) influences the choice of the $\varepsilon(i,j)$. Below, we formulate a criterion for the polynomials $\varepsilon(i,j)$ that leads to an efficient implementation of SoS. However, we do not claim that other choices of the $\varepsilon(i,j)$ cannot lead to efficient implementations too.

Recall that a primitive operation is a function $f$ that maps a set $Q$ of $k$ objects to $+1$, 0, or $-1$. If the $\varepsilon$-expansion is defined properly, then $f(Q(\varepsilon)) \in \{+1,-1\}$ provided $\varepsilon > 0$ is small enough. In general, $f(Q(\varepsilon))$ will be the sign of a fairly complicated function in $\varepsilon$. (Since $f$ is now a binary function we can identify $\{+1,-1\}$ with $\{\text{true}, \text{false}\}$ and express it as a predicate. We will follow this practice in the following sections of this paper.) One way to allow for an efficient evaluation of $f(Q(\varepsilon))$ is to choose the $\varepsilon(i,j)$ in different orders of magnitude such that two expressions, each consisting of several factors of the form $\varepsilon(i,j)$, can be compared solely on the basis of the index pairs $(i,j)$ involved. When we evaluate $f(Q(\varepsilon))$, we can sort its terms in order of decreasing significance which can be done by comparing sets of index pairs. The most significant term will be a term without any $\varepsilon$-factor; it will be equal to $f(Q)$. The first term with a nonzero coefficient decides the sign of the function. If $Q$ is nondegenerate to begin with, then $f(Q(\varepsilon)) = f(Q)$, and no other term has to be determined. In Sections 3 and 4, we will see that such a choice of the $\varepsilon(i,j)$ allows us to determine the sign of a fairly complicated polynomial in only a few steps.

Note that SoS requires us to tell when $Q$ is degenerate, which means that we need to be able to decide whether or not $f(Q) = 0$. This is not possible with the kind of floating-point arithmetic that it is usually provided by current computers. Instead, we need to use exact arithmetic and, thus, occasionally long integers. These admittedly somewhat expensive operations occur only inside the primitives and do not concern the user of SoS. Furthermore, the length of such long integers is bounded by a constant if $kd$, the number of input parameters of $f$, is bounded by a constant. In most geometric algorithms, this constant is reasonably small. In Section 6 we report on our experience in implementing SoS and give an indication to what extent the use of long integer arithmetic slows down the computation. This point cannot be taken lightly because the long integer arithmetic is likely to occur in the innermost loop of any program that uses SoS and thus dictates the constant in front of the asymptotic running time. However, it is worthwhile to mention that the need for exact arithmetic is not a peculiar feature of SoS itself, but is necessary whenever we do exact computation rather than push our luck and hope for the cancellation of round-off errors.

# 3 Finite Point Sets — a Case Study

For a further discussion of SoS it is advantageous to apply it to certain geometric objects and certain primitive operations defined for these objects. We choose points in the $d$-dimensional Euclidean space $E^d$ as the objects for the case study. Notice that this is actually no loss of generality since every object specified by $d$ parameters can be interpreted as a point in $E^d$. The primitive operation that we will consider takes $d+1$ points as input and decides on which side of the hyperplane spanned by the last $d$ points the first point lies. As we will see in Section 5, this primitive operation has a wide range of applications.

If a given finite point set is perturbed, as explained in Section 2, one can ignore all degeneracies and special cases. The price for this simulated simplicity is that the coordinates of the points are now symbolic expressions in $\varepsilon$. Even for a simple task, such as the comparison of two coordinates, we need a custom-made procedure that handles the $\varepsilon$-expansions of the coordinates. Let $\pi_{i,j}$ be the $j$-th coordinate of point $p_i$ and let $\pi_{k,l}$ be the $l$-th coordinate of $p_k$, $0 \le i,k \le n-1$ and $1 \leq j, l \leq d$. To decide which one of the two corresponding perturbed coordinates is smaller, we define a predicate *Smaller* as follows:

$$\mathit{Smaller}(\pi_{i,j}; \pi_{k,l}) = \mathit{true} \quad \text{iff} \quad \pi_{i,j}(\varepsilon) < \pi_{k,l}(\varepsilon).$$

Due to SoS, we can neglect degeneracies, i.e., we have $\pi_{i,j}(\varepsilon) \neq \pi_{k,l}(\varepsilon)$, and for this reason the predicate $\mathit{Smaller}(\pi_{i,j}; \pi_{k,l}) = \mathit{false}$ if and only if $\pi_{i,j}(\varepsilon) > \pi_{k,l}(\varepsilon)$. The implementation of this predicate is fairly straightforward since we can compare the $\varepsilon$-terms, $\varepsilon(i, j)$ and $\varepsilon(k, l)$, by comparing the defining index pairs (see Section 3.2, Lemma 3.2).

**Predicate 1** (*Smaller*) Assume the $\varepsilon$-expansion $\varepsilon(i, j)$ defined as in Section 3.2 (2). With this, for indices $0 \leq i, k \leq n-1$ and $1 \leq j, l \leq d$ which satisfy $(i, j) \neq (k, l)$, the predicate $\mathit{Smaller}(\pi_{i,j}; \pi_{k,l})$ can be implemented as follows.

**function** *Smaller* $(\pi_{i,j}; \pi_{k,l})$ **returns** Boolean
**begin**
  **if** $\pi_{i,j} \neq \pi_{k,l}$ **then**
    **return** $(\pi_{i,j} < \pi_{k,l})$
  **else if** $i \neq k$ **then**
    **return** $(i > k)$
  **else**
    **return** $(j < l)$
**end**

Notice that, in this case, the coordinates $\pi_{i,j}$ and $\pi_{k,l}$ as well as their index pairs $(i, j)$ and $(k, l)$ have to be passed as arguments whenever predicate *Smaller* is called. This means that in popular programming languages, such as *Pascal*, the function heading would be something like

```
FUNCTION smaller (i, j, k, l, Pij, Pkl): Boolean;
```

but implementation details like this will be ignored in the remainder of the paper. Furthermore, notice that we have

$$\mathit{Smaller}(\pi_{i,j}; \pi_{k,l}) = \mathit{true} \quad \text{iff} \quad \det \begin{pmatrix} \pi_{i,j}(\varepsilon) & 1 \\ \pi_{k,l}(\varepsilon) & 1 \end{pmatrix} < 0.$$

In Section 3.1, we will express more complicated predicates than just comparisons of coordinates by similar determinants. For matrices not exceeding a given size it is not difficult to specify the $\varepsilon$-expansion $\varepsilon(i, j)$ such that all requirements discussed in Section 2 are satisfied. This will be done in Section 3.2. Finally, Section 3.3 extends the results to homogeneous coordinates. The procedures that implement the predicates will be developed in Section 4.

## 3.1 Predicates Expressed by Determinants

This section introduces the notion of orientation of a sequence of $d + 1$ points in $E^d$. With this concept, we will be able to give an implementation of the primitive operation for $d + 1$ points mentioned above.

The *orientation* of a sequence of points $(p_{i_0}, p_{i_1}, \ldots, p_{i_d})$ in $E^d$ is either *negative* or *positive* — unless the $d+1$ points lie in a common hyperplane, in which case the orientation is undefined. The exceptional case is a degeneracy that can be ignored if the points are perturbed. We define the orientation of a sequence recursively. It will be important that the orientation of a sequence depends only on the relative position of the points to each other and not on their absolute positions.

If the dimension $d = 1$, then the orientation of $(p_{i_0}, p_{i_1})$ is positive if $p_{i_0} > p_{i_1}$ and it is negative if $p_{i_0} < p_{i_1}$ (compare with Figure 3-I, (a) and (b)). If $d = 2$, then $(p_{i_0}, p_{i_1}, p_{i_2})$ has positive orientation if the three points define a left-turn in the plane, that is, $p_{i_2}$ lies to the left of the directed line that passes through $p_{i_0}$ and $p_{i_1}$ in this order. If $(p_{i_0}, p_{i_1}, p_{i_2})$ defines a right-turn, then its orientation is negative. Note that the orientation of $(p_{i_0}, p_{i_1}, p_{i_2})$ is the same as the orientation of $(p_{i_1}, p_{i_2})$ as "seen from" $p_{i_0}$. Indeed, the line through $p_{i_1}$ and $p_{i_2}$ can be identified with $E^1$ as soon as we choose a direction of the line. This direction is provided by the location of $p_{i_0}$: It goes from left to right as seen from $p_{i_0}$ (see Figure 3-I, (c) and (d)).

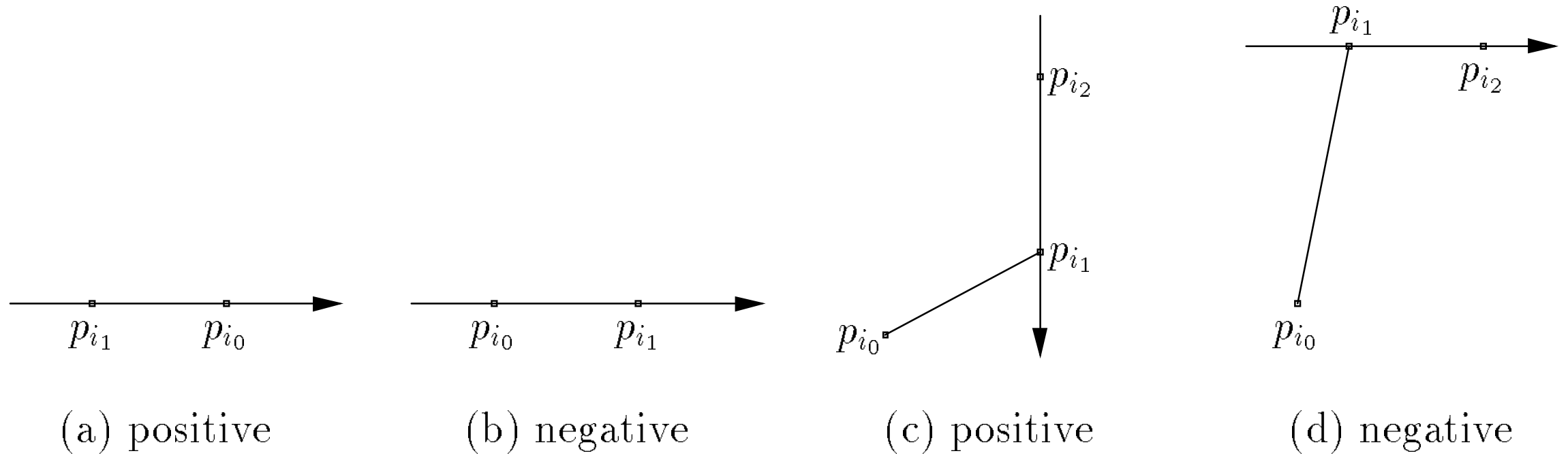

Figure 3-I: The orientation of $d+1$ points in dimension $d$, for $d = 1, 2$.

If $d > 2$, then the orientation of $(p_{i_0}, p_{i_1}, \ldots, p_{i_d})$ is the same as the orientation of $(p_{i_1}, \ldots, p_{i_d})$ as seen from $p_{i_0}$. For example, $(p_{i_0}, p_{i_1}, p_{i_2}, p_{i_3})$ in $E^3$ has positive orientation if $p_{i_0}$ observes $(p_{i_1}, p_{i_2}, p_{i_3})$ making a left-turn. In most situations where the concept of orientation is used, the interest is in the position of one point, $p_{i_0}$, relative to $d$ other points, $p_{i_1}, p_{i_2}, \ldots, p_{i_d}$. We thus say that $p_{i_0}$ lies on the *positive side of* $(p_{i_1}, \ldots, p_{i_d})$ if $(p_{i_0}, p_{i_1}, \ldots, p_{i_d})$ has positive orientation, and $p_{i_0}$ lies on the *negative side of* $(p_{i_1}, \ldots, p_{i_d})$ if $(p_{i_0}, p_{i_1}, \ldots, p_{i_d})$ has negative orientation.

To decide upon the orientation of a sequence of $d+1$ points in $E^d$, we use the matrix

$$\Lambda = \begin{pmatrix} \pi_{i_0,1} & \pi_{i_0,2} & \cdots & \pi_{i_0,d} & 1 \\ \pi_{i_1,1} & \pi_{i_1,2} & \cdots & \pi_{i_1,d} & 1 \\ \vdots & \vdots & \ddots & \vdots & \vdots \\ \pi_{i_d,1} & \pi_{i_d,2} & \cdots & \pi_{i_d,d} & 1 \end{pmatrix}. \tag{3-a}$$

**Lemma 3.1** The orientation of $(p_{i_0}, p_{i_1}, \ldots, p_{i_d})$ is positive if and only if $\mathrm{sign}(\det \Lambda) = +1$ and it is negative if and only if $\mathrm{sign}(\det \Lambda) = -1$.

Notice that $\det \Lambda$ vanishes if and only if the $d+1$ points are degenerate, that is, they lie in a common hyperplane — a case that can be neglected within the perturbed point set $P(\varepsilon)$. Recall from linear algebra that the determinant of a matrix is multiplied by $-1$ if we exchange two rows. Thus, the orientation of a permutation of $(p_{i_0}, p_{i_1}, \ldots, p_{i_d})$ is the same as the orientation of the sequence itself if the number of transpositions is even, otherwise, its orientation is the opposite of the orientation of $(p_{i_0}, p_{i_1}, \ldots, p_{i_d})$.

There are plenty of algorithms for point set problems which are based on computing the orientation of a sequence of points. Prime examples are the construction of convex hulls (see [PH77], [PS85], [Se81], [Se86], or [Ed87]), computing $\lambda$-matrices as discussed in [GP83] and [Ed87], and finding convex subsets (see [CK80], [EG89], and [Ed87]). The remainder of this section considers the primitive operations required by the three-dimensional convex hull algorithm of Preparata and Hong which is described in [PH77], [PS85], and [Ed87].

The first step of the algorithm sorts the points in $x_1$-direction. To perform this step, it needs to compare the $x_1$-coordinates of two points, which can be done by computing the orientation of their orthogonal projections onto the $x_1$-axis. Second, it constructs the two-dimensional convex hull of the points projected onto the $x_1x_3$-plane. Here, the primitive operation is to decide whether three points (in the $x_1x_3$-plane) define a left-turn or a right-turn. Third, the algorithm constructs the three-dimensional convex hull by repeating the following operation:

> Given a plane pivoting about two extreme points $p_{i_3}$ and $p_{i_2}$, find the point hit first by this plane.

This operation can be reduced to a number of comparisons of the form: Given two points $p_{i_1}$ and $p_{i_0}$, which one is hit earlier by the pivoting plane? To perform such a comparison is equivalent to deciding on which side of the plane through $p_{i_3}$, $p_{i_2}$, and $p_{i_1}$ point $p_{i_0}$ lies. This is the same as computing the orientation of $(p_{i_0}, p_{i_1}, p_{i_2}, p_{i_3})$. Thus, we see that the convex hull algorithm of Preparata and Hong requires three primitive operations all of which determine the orientation of point sequences.

## 3.2 Choosing the Form of the Perturbation

As explained in Section 3.1 the primitive operation that determines the orientation of a sequence of $d+1$ points in $d$ dimensions computes the sign of a determinant of a $(d+1)$-by-$(d+1)$ matrix. SoS replaces the coordinates $\pi_{i,j}$ in this matrix by entries of the form $\pi_{i,j}+\varepsilon(i,j)$. The determinant itself is then the sum of a finite number of terms, where each term is the product of $d$ items and an item is either an original coordinate or an $\varepsilon(i,j)$. Thus, each term consists of a coefficient, which is the product of original coordinates, and a so-called *$\varepsilon$-product*, a product of factors of the form $\varepsilon(i,j)$. The number of factors $\varepsilon(i,j)$ can be zero in which case the $\varepsilon$-product is defined to be equal to 1. As mentioned in Section 2, it is irrelevant what exactly the definition of the $\varepsilon$-expansion is as long as it satisfies certain requirements. The computational simulation is uneffected if we change the definition of the $\varepsilon$-expansion within allowed limits. Even so, it is important to show that there is at least one $\varepsilon$-expansion that satisfies the requirements. The existence of such an expansion implies the physical existence of an appropriately perturbed point set which is the only guarantee for the consistency of our method we have.

We define

$$\varepsilon(i,j) = \varepsilon^{2^{i\cdot\delta-j}}, \tag{3-b}$$

for $0 \le i \le n-1$, $1 \le j \le d$, and $\delta \ge d$, and show that this choice satisfies all the requirements of SoS. Notice that the amount of perturbation experienced by coordinate $\pi_{i,j}$ is larger than the perturbation of $\pi_{k,l}$ if and only if $(i,j) \prec (k,l)$, that is, $i < k$ or $i = k$ and $j > l$. Furthermore, we have

$$\prod_{(i,j)\prec(k,l)} \varepsilon(i,j) = \prod_{(i,j)\prec(k,l)} \varepsilon^{2^{i\cdot\delta-j}} > \varepsilon^{2^{k\cdot\delta-l}} = \varepsilon(k,l) \tag{3-c}$$

if $0 < \varepsilon < 1$. This is equivalent to stating that $2^{k \cdot \delta - l}$, the exponent of $\varepsilon(k,l)$, is larger than the sum of the exponents of all $\varepsilon(i,j)$ with $(i,j) \prec (k,l)$. It follows that it is sufficient to consider the sets of index pairs when we compare two $\varepsilon$-products. Let $e_1$ and $e_2$ be two different $\varepsilon$-products and let $\mathcal{I}(e_1)$ and $\mathcal{I}(e_2)$ be the two associated sets of index pairs. We call $\mathcal{I}(e_1)$ *smaller* than $\mathcal{I}(e_2)$ if the set $\mathcal{I}(e_1) - \mathcal{I}(e_2)$ is empty or if $(i,j) \prec (k,l)$, for $(i,j)$ the largest index pair in $\mathcal{I}(e_1) - \mathcal{I}(e_2)$ and $(k,l)$ the largest index pair in $\mathcal{I}(e_2) - \mathcal{I}(e_1)$.

**Lemma 3.2** Let $c_1$ and $c_2$ be two positive constants and let $e_1$ and $e_2$ be two different $\varepsilon$-products. Then $c_1 \cdot e_1 > c_2 \cdot e_2$ for a small enough $\varepsilon$ if $\mathcal{I}(e_1)$ is smaller than $\mathcal{I}(e_2)$.

Lemma 3.2 is an immediate consequence of (3-c) and the fact that a small enough $\varepsilon$ can compensate the influence of the constants $c_1$ and $c_2$. Notice that it is actually irrelevant which index pairs $\mathcal{I}(e_1)$ and $\mathcal{I}(e_2)$ contain. The only thing of importance is the relative position of $\mathcal{I}(e_1)$ and $\mathcal{I}(e_2)$ in the ordering of all sets of index pairs, where large index pairs are more significant in the comparison of sets than small index pairs. Observe also that Lemma 3.2 holds if we increase the value of $\delta$ in the definition of the $\varepsilon$-expansion. It turns out that this lemma is the crucial property that allows us to prove that $P(\varepsilon)$, the perturbed point set, is simple and that the orientation of $d+1$ points in $P(\varepsilon)$ can be computed efficiently.

**Lemma 3.3** The set $P(\varepsilon)$ is nondegenerate if $\varepsilon > 0$ is sufficiently small.

**Proof.** To prove the assertion, we show that for no choice of $d+1$ mutually distinct indices $i_0, i_1, \ldots, i_d$, the determinant of the matrix

$$\Lambda(\varepsilon) = \begin{pmatrix} \pi_{i_0,1} + \varepsilon^{2^{i_0 \cdot \delta - 1}} & \pi_{i_0,2} + \varepsilon^{2^{i_0 \cdot \delta - 2}} & \cdots & \pi_{i_0,d} + \varepsilon^{2^{i_0 \cdot \delta - d}} & 1 \\ \pi_{i_1,1} + \varepsilon^{2^{i_1 \cdot \delta - 1}} & \pi_{i_1,2} + \varepsilon^{2^{i_1 \cdot \delta - 2}} & \cdots & \pi_{i_1,d} + \varepsilon^{2^{i_1 \cdot \delta - d}} & 1 \\ \vdots & \vdots & \ddots & \vdots & \vdots \\ \pi_{i_d,1} + \varepsilon^{2^{i_d \cdot \delta - 1}} & \pi_{i_d,2} + \varepsilon^{2^{i_d \cdot \delta - 2}} & \cdots & \pi_{i_d,d} + \varepsilon^{2^{i_d \cdot \delta - d}} & 1 \end{pmatrix} \tag{3-d}$$

is equal to zero. To see this, we assume w.l.o.g. that $0 \leq i_0 < i_1 < \ldots < i_d \leq n-1$ and sort the terms of $\det \Lambda(\varepsilon)$ in order of increasing exponents of $\varepsilon$. Specifically,

$$\det \Lambda$$

is the first term, and

$$(-1)^{\lceil d/2 \rceil} \cdot \varepsilon^{2^{i_1 \cdot \delta - d} + 2^{i_2 \cdot \delta - (d-1)} + \ldots + 2^{i_d \cdot \delta - 1}}$$

the last one. Each term is of the form $b \cdot \varepsilon^c$, for some constants $b$ and $c$. Because we can assume that $\varepsilon > 0$ is arbitrarily small, the absolute value of the first term with nonzero coefficient $b$ is bigger than the sum of all other terms. Furthermore, such a term always exists since (3-c) guarantees that no two terms of the determinant have the same exponent of $\varepsilon$ and thus such a term cannot cancel. For example, the coefficient of the last term is $(-1)^{\lceil d/2 \rceil} \neq 0$ and cannot be canceled by any other term. Consequently, $\det \Lambda(\varepsilon)$ does not vanish. □

As pointed out in the proof of Lemma 3.3, the most significant term of the polynomial $\det \Lambda(\varepsilon)$ is the determinant $\det \Lambda$ of the original coordinates. If the orientation of the original sequence $(p_{i_0}, p_{i_1}, \ldots, p_{i_d})$ is defined, then this term is nonzero which implies that the orientation of the perturbed sequence is the same. This is reassuring since it shows that the perturbation does not change nondegenerate relations of the original point set.

The curious reader might wonder why the perturbation is defined in the peculiar form given by the $\varepsilon$-expansion (3-b). As mentioned before, there are many other choices that could be used, for example,

$$\varepsilon(i,j) = \varepsilon^{2^{i \cdot \delta + j}}$$

is such a possibility. This $\varepsilon$-expansion would also work but its implementation is slightly more difficult than that of (3-b) (compare with Section 4.2). On the other hand, many less "exotic" choices do not work. The remainder of this section illustrates this by considering two choices of $\varepsilon(i,j)$ which appear simpler than (3-b). The two choices are

$$\varepsilon(i,j) = \varepsilon^{i \cdot \delta + j} \quad \text{and} \quad \varepsilon(i,j) = (i \cdot \delta + j) \cdot \varepsilon.$$

In both cases, Lemma 3.3 does not hold. The reason for the failure is that both expansions do not satisfy (3-c) and thus possibly lead to cancellations of $\varepsilon$-terms in $\det \Lambda(\varepsilon)$. Such cancellations occur for example if all $d+1$ points of the sequence coincide with the origin. In this case the matrix $\Lambda(\varepsilon)$ equals

$$\begin{pmatrix} \varepsilon(i_0,1) & \varepsilon(i_0,2) & \cdots & \varepsilon(i_0,d) & 1 \\ \varepsilon(i_1,1) & \varepsilon(i_1,2) & \cdots & \varepsilon(i_1,d) & 1 \\ \vdots & \vdots & \ddots & \vdots & \vdots \\ \varepsilon(i_d,1) & \varepsilon(i_d,2) & \cdots & \varepsilon(i_d,d) & 1 \end{pmatrix}.$$

If we define $\varepsilon(i,j) = \varepsilon^{i \cdot \delta + j}$, then the second column is equal to $\varepsilon$ times the first column which implies that $\det \Lambda(\varepsilon) = 0$ if $d \geq 2$. If $\varepsilon(i,j) = (i \cdot \delta + j) \cdot \varepsilon$, then the sum of the first and the third columns equals twice the second column; hence, $\det \Lambda(\varepsilon) = 0$ if $d \geq 3$.

## 3.3 Homogeneous Coordinates

When we develop the primitive procedures for computing the orientation of $d+1$ points in Section 4, we represent a point by its homogeneous coordinates. This representation is slightly more general than ordinary Cartesian coordinates (it can also represent points at infinity) and leads to a slightly more uniform procedural treatment.

Let $p$ be a point in $E^d$ and let $(\pi_1^C, \pi_2^C, \ldots, \pi_d^C)$ be its sequence of Cartesian coordinates. Point $p$ has $d+1$ *homogeneous coordinates*

$$(\pi_1^H, \pi_2^H, \ldots, \pi_d^H; \pi_{d+1}^H)$$

such that

$$\pi_i^C = \frac{\pi_i^H}{\pi_{d+1}^H}, \quad \text{for} \quad 1 \leq i \leq d.$$

Thus, $p$ is $1/\pi_{d+1}^H$ times the point whose Cartesian coordinates are equal to the first $d$ homogeneous coordinates of $p$. Notice that the homogeneous coordinates of $p$ are not unique; we still represent the same point $p$ if we multiply each coordinate by the same nonzero scalar. If we decrease the absolute value of $\pi_{d+1}^H$ without changing the other homogeneous coordinates, then $p$ moves away from the origin on a straight line and it reaches "infinity" when $\pi_{d+1}^H$ becomes 0. Indeed, $p$ is "at infinity" if and only if $\pi_{d+1}^H = 0$. Using homogeneous coordinates, it is not allowed to have all $d+1$ coordinates are equal to 0 — in this event $p$ is not defined.

We next extend Lemma 3.1 to homogeneous coordinates, that is, we characterize the orientation of a sequence of $d+1$ points $(p_{i_0}, p_{i_1}, \ldots, p_{i_d})$,

$$p_{i_\nu} = (\pi^H_{i_\nu,1}, \pi^H_{i_\nu,2}, \ldots, \pi^H_{i_\nu,d}; \pi^H_{i_\nu,d+1})$$

in terms of their homogeneous coordinates. The orientation of a sequence of $d+1$ points is not defined if any of the points lies at infinity. In fact, it is not possible to generalize the notion of orientation to points at infinity without changing our interpretation of a point at infinity. For example, consider a sequence $S$ of $d$ finite points and one point $p = (\pi^H_1, \pi^H_2, \ldots, \pi^H_d; 0)$ at infinity. We can think of $p$ as the limit of points

$$p(\epsilon) = (\pi^H_1, \pi^H_2, \ldots, \pi^H_d; \epsilon),$$

when $\epsilon > 0$ goes to zero, but as well, we can think of $p$ as the limit of these points if $\epsilon$ is negative and approaches zero. If we replace $p$ by $p(\epsilon)$ with $\epsilon$ small enough, then $\epsilon > 0$ and $\epsilon < 0$ lead to different orientations. We thus restrict our discussion of orientation to finite points. Define

$$\Delta = \begin{pmatrix} \pi^H_{i_0,1} & \pi^H_{i_0,2} & \cdots & \pi^H_{i_0,d+1} \\ \pi^H_{i_1,1} & \pi^H_{i_1,2} & \cdots & \pi^H_{i_1,d+1} \\ \vdots & \vdots & \ddots & \vdots \\ \pi^H_{i_d,1} & \pi^H_{i_d,2} & \cdots & \pi^H_{i_d,d+1} \end{pmatrix}. \tag{3-e}$$

If $\pi^H_{i_\nu,d+1} = 1$, for $0 \leq \nu \leq d$, then $\Delta$ is the same as the matrix $\Lambda$ used in Lemma 3.1. Otherwise, we can multiply the rows such that $\pi^H_{i_\nu,d+1} = 1$. The sign of $\det \Delta$ changes if we multiply a row with a negative number, which implies the following result:

**Lemma 3.4** Let $(p_{i_0}, p_{i_1}, \ldots, p_{i_d})$ be a sequence of points with $p_{i_\nu} = (\pi^H_{i_\nu,1}, \pi^H_{i_\nu,2} \ldots, \pi^H_{i_\nu,d}, \pi^H_{i_\nu,d+1})$ and $\pi^H_{i_\nu,d+1} \neq 0$. Their orientation is positive if $\text{sign}(\det \Delta) = \prod_{\nu=0}^{d} \text{sign}(\pi^H_{i_\nu,d+1})$, negative if $\text{sign}(\det \Delta) = -\prod_{\nu=0}^{d} \text{sign}(\pi^H_{i_\nu,d+1})$, and undefined if $\det \Delta = 0$.

In contrast to Cartesian coordinates, a point is now represented by $d+1$ coordinates which makes it necessary to choose $\delta \geq d+1$ when defining the $\varepsilon$-expansion $\varepsilon(i,j)$ in (3-b). With this, it is easy to prove that determinants cannot vanish which implies that Lemma 3.3 holds also for the new setting using homogeneous coordinates.

# 4 Implementing a Predicate

This section presents the actual implementation of a geometric predicate using SoS. The chosen predicate determines the orientation of a sequence of points, as defined in Section 3. Its implementation will be based on the $\varepsilon$-expansion specified in Section 3.2 (3-b) and on the fact that the orientation can be found by evaluating the sign of a determinant as stated in Sections 3.1 and 3.3. The crux of the implementation is that this determinant is a polynomial in $\varepsilon$. The computation of the sign of such a polynomial is discussed in Section 4.1. The coefficients of the polynomial turn out to be subdeterminants of the original matrix. Based on this observation, Section 4.2 gives an algorithm that generates these subdeterminants in sequence of decreasing significance by employing a special encoding scheme. Finally, in Section 4.3 we will briefly address the problem of sign computation of integer determinants in general.

In Sections 3.1 and 3.3 we defined the "orientation" of a sequence of points in $d$-dimensional Euclidean space given by Cartesian and homogeneous coordinates. We now formally develop the corresponding predicate that uses perturbation in the sense of SoS. In the Cartesian case each point $p_\nu$ is given by its $d$ coordinates

$$p_\nu = (\pi_{\nu,1}, \ldots, \pi_{\nu,d}),$$

whereas in the case of homogeneous coordinates a point is represented by a $(d+1)$-tupel

$$p_\nu = (\pi_{\nu,1}, \ldots, \pi_{\nu,d}; \pi_{\nu,d+1}).$$

Let

$$P = \{p_0, \ldots, p_{n-1}\}$$

be a set of $n$ points in $E^d$, and denote by

$$P(\varepsilon) = \{p_0(\varepsilon), \ldots, p_{n-1}(\varepsilon)\}$$

its perturbed version using the $\varepsilon$-expansion of Section 3.2 (3-b), assuming $\delta$ large enough such that Lemma 3.2 is valid. Now define for $d+1$ points with distinct indices $i_0, i_1, \ldots, i_d$, all in the range from 0 through $n-1$,

$$\mathit{Positive}_d(p_{i_0}, \ldots, p_{i_d}) = \text{true} \quad \text{iff} \quad \text{the orientation of } (p_{i_0}(\varepsilon), \ldots, p_{i_d}(\varepsilon)) \text{ is positive.}$$

Degenerate cases can be neglected because we simulate simplicity. From Lemma 3.1 it follows that $\mathit{Positive}_d$ is equivalent to the test whether or not

$$\text{sign}(\det \Lambda(\varepsilon)) = +1,$$

with $\Lambda(\varepsilon)$ denoting the the corresponding matrix of the perturbed Cartesian coordinates as in (3-d). In the homogeneous case (see Lemma 3.4) we have to check whether or not

$$\text{sign}(\det \Delta(\varepsilon)) = \prod_{\nu=0}^{d} \text{sign}(\pi_{i_\nu,d+1}(\varepsilon)).$$

Here, $\Delta(\varepsilon)$ denotes the perturbed version of matrix $\Delta$ in (3-e), whose rows are formed by the homogeneous coordinates of the points involved; that is,

$$\Delta(\varepsilon) = \begin{pmatrix} \pi_{i_0,1} + \varepsilon(i_0,1) & \pi_{i_0,2} + \varepsilon(i_0,2) & \cdots & \pi_{i_0,d+1} + \varepsilon(i_0,d+1) \\ \pi_{i_1,1} + \varepsilon(i_1,1) & \pi_{i_1,2} + \varepsilon(i_1,2) & \cdots & \pi_{i_1,d+1} + \varepsilon(i_1,d+1) \\ \vdots & \vdots & \ddots & \vdots \\ \pi_{i_d,1} + \varepsilon(i_d,1) & \pi_{i_d,2} + \varepsilon(i_d,2) & \cdots & \pi_{i_d,d+1} + \varepsilon(i_d,d+1) \end{pmatrix}.$$

At first sight, the development of such an $\varepsilon$-determinant seems to be a painful exercise. Yet, it will turn out that it is not that hard and can be achieved in an algorithmically clean way. Anyway, to begin with something easy, consider

$$\det \Delta_2(\varepsilon) = \det \begin{pmatrix} \pi_{i,1} + \varepsilon(i,1) & \pi_{i,2} + \varepsilon(i,2) \\ \pi_{j,1} + \varepsilon(j,1) & \pi_{j,2} + \varepsilon(j,2) \end{pmatrix}.$$

Let $\varepsilon((i_1, j_1), \ldots, (i_k, j_k)) = \prod_{\nu=1}^{k} \varepsilon(i_\nu, j_\nu)$, and call it a *$k$-fold $\varepsilon$-product*; $\varepsilon() = 1$ is called the *0-fold $\varepsilon$-product*. Furthermore assume $i < j$. When we now develop the determinant we get

$$\begin{aligned} \det \Delta_2(\varepsilon) = \; & + \begin{pmatrix} \pi_{i,1} & \pi_{i,2} \\ \pi_{j,1} & \pi_{j,2} \end{pmatrix} \cdot \varepsilon() - \\ & - \pi_{j,1} \cdot \varepsilon(i, 2) + \pi_{j,2} \cdot \varepsilon(i, 1) + \\ & + \pi_{i,1} \cdot \varepsilon(j, 2) + 1 \cdot \varepsilon((j, 2), (i, 1)) - \\ & - \pi_{i,2} \cdot \varepsilon(j, 1) - 1 \cdot \varepsilon((j, 1), (i, 2)), \end{aligned} \tag{4-a}$$

where the terms are already sorted by increasing powers of $\varepsilon$. Note again that the first coefficient corresponds to the "unperturbed" determinant, i.e., $\Delta_2$, whose evaluation would be part of any implementation of the predicate — of course, followed by the more or less awkward handling of all possible degeneracies. Observe also that the coefficient of the fifth term is a constant, namely $+1$. Thus, the last two terms have no influence on the sign of $\det \Delta_2(\varepsilon)$. Therefore, the number of relevant terms of the $\varepsilon$-polynomial $\det \Delta_2(\varepsilon)$ is only 5, rather than 7 which is the total number of terms.

It is convenient to assume $i_0 < \ldots < i_d$ (compare with (4-a)). This assumption together with Lemma 3.2 implies that the sign of $\det \Delta(\varepsilon)$ and $\det \Lambda(\varepsilon)$ can be computed without any further knowledge of the values of the indices. Clearly, this is not the case in general but can always be achieved by appropriate row exchanges in $\Delta(\varepsilon)$ or $\Lambda(\varepsilon)$ — recall that each exchange changes the sign of the determinant. For this, assume a procedure $\mathit{Sort}_{d+1}((i_0, \ldots, i_d), (i'_0, \ldots, i'_d), s')$ that returns for a given sequence of $d+1$ indices $(i_0, \ldots, i_d)$ the sorted sequence $(i'_0, \ldots, i'_d)$. Additionally, $\mathit{Sort}_{d+1}$ returns $s'$ which is set to the number of exchanges used. We can now implement predicate $\mathit{Positive}_d$ using two operations, $\mathit{SignDet}\Lambda$ and $\mathit{SignDet}\Delta$, that compute the sign of the $\varepsilon$-polynomials $\det \Lambda(\varepsilon)$ and $\det \Delta(\varepsilon)$ assuming $i_0 < \ldots < i_d$. Both functions will be discussed in Section 4.1.

**Predicate 2** (*Positive*) Let $p_{i_0}, \ldots, p_{i_d}$ be $d+1$ points in $E^d$ given in Cartesian or homogeneous coordinates with distinct indices all between 0 and $n-1$. Then the following pseudocode is an implementation of the predicate $\mathit{Positive}_d$:

**function** $Positive_d\ (p_{i_0}, \ldots, p_{i_d})$ **returns** Boolean
**local** $i'_0, \ldots, i'_d$, $d'$, $s'$, $\nu$
**begin**
$Sort_{d+1}((i_0, \ldots, i_d), (i'_1, \ldots, i'_d), s')$
**if** Cartesian coordinates **then**

$$d' \leftarrow SignDet\Lambda_{d+1}\begin{pmatrix} \pi_{i'_0,1}(\varepsilon) & \cdots & \pi_{i'_0,d}(\varepsilon) & 1 \\ \vdots & \ddots & \vdots & \vdots \\ \pi_{i'_d,1}(\varepsilon) & \cdots & \pi_{i'_d,d}(\varepsilon) & 1 \end{pmatrix}$$

**else**

$$d' \leftarrow SignDet\Delta_{d+1}\begin{pmatrix} \pi_{i'_0,1}(\varepsilon) & \cdots & \pi_{i'_0,d}(\varepsilon) & \pi_{i'_0,d+1}(\varepsilon) \\ \vdots & \ddots & \vdots & \vdots \\ \pi_{i'_d,1}(\varepsilon) & \cdots & \pi_{i'_d,d}(\varepsilon) & \pi_{i'_d,d+1}(\varepsilon) \end{pmatrix}$$

**if** odd($s'$) **then** $d' \leftarrow -d'$
**if** Cartesian coordinates **then**
**return** $(d' = +1)$
**else**
**return** $(d' = \prod_{\nu=0}^{d} \operatorname{sign}(\pi_{i_\nu,d+1}(\varepsilon)))$
**end**

The problem is now to give efficient implementations for the two functions $SignDet\Lambda_{d+1}$ and $SignDet\Delta_{d+1}$. We feel that it is important to stress that "efficiency" is meant in a practical sense — in theory it can be done in constant time anyway, assuming $d$ is a constant.

## 4.1 The Sign of a Perturbed Determinant

We now illustrate the implementation of SoS on the bottommost programming level by implementing function $SignDet\Delta_D$, which returns the sign of a $D$-by-$D$ $\varepsilon$-determinant $\det \Delta_D(\varepsilon)$ for any given $D$; primitive $SignDet\Lambda_D$ can be treated in the same way. To appreciate the significance of a (practically) efficient implementation of $SignDet\Delta_D$ we point out that this is in fact the major part of SoS, at least when applied to the predicate described above. Provided that $i_0 < \cdots < i_D$, we will show that it is possible without great effort to generate the sequence of the coefficients of $\det \Delta_D(\varepsilon)$ in decreasing order of significance. Since $\varepsilon$ can be assumed to be sufficiently small (but positive), the sign of the $\varepsilon$-polynomial is therefore equivalent to the sign of the first nonvanishing coefficient.

Using simple rules for evaluating a determinant as exemplified for $\det \Delta_2(\varepsilon)$ in (4-a), the coefficient of every term in $\det \Delta_D(\varepsilon)$ is a subdeterminant of the "unperturbed" matrix $\Delta_D$. Here, a single entry is called a 1-by-1 subdeterminant and, by definition, the 0-by-0 subdeterminant is equal to 1. To tell the whole truth, we must mention that each coefficient in effect is a subdeterminant together with a certain sign, that is, multiplied by either $+1$ or $-1$. We will see in Section 4.2 how to decide whether $+1$ or $-1$ applies. To continue our discussion, we need a few notations. We say that the $(t+1)$-st coefficient in order of decreasing significance, denoted by $\det M_t^{\Delta_D}$, is the *cofactor of depth* $t$ of matrix $\Delta_D(\varepsilon)$. Note that this coefficient already includes its proper sign. Thus, $\det M_0^{\Delta_D} = +\det \Delta_D$. The size of the corresponding matrix (i.e., the number of rows or columns) is denoted by $k_t = k(M_t^{\Delta_D})$. These definitions are illustrated in Table 1 which shows all significant terms of $\det \Delta_2(\varepsilon)$. In the column with heading $\varepsilon_t$ we display the $\varepsilon$-product associated with the cofactor of depth $t$. The column $v_t$ will be explained later.

| $t$ | $k_t \cdot k_t$ | $v_t$ | $\det M_t^{\Delta_2}$ | $\varepsilon_t$ |
|---|---|---|---|---|
| 0 | $2 \cdot 2$ | $[3,3;3]$ | $+\det \begin{pmatrix} \pi_{i,1} & \pi_{i,2} \\ \pi_{j,1} & \pi_{j,2} \end{pmatrix}$ | $\varepsilon()$ |
| 1 | $1 \cdot 1$ | $[2,3;3]$ | $-\det(\pi_{j,1}) = -\pi_{j,1}$ | $\varepsilon(i,2)$ |
| 2 | $1 \cdot 1$ | $[1,3;3]$ | $+\det(\pi_{j,2}) = +\pi_{j,2}$ | $\varepsilon(i,1)$ |
| 3 | $1 \cdot 1$ | $[2,2;3]$ | $+\det(\pi_{i,1}) = +\pi_{i,1}$ | $\varepsilon(j,2)$ |
| 4 | $0 \cdot 0$ | $[1,2;3]$ | $+\det() = +1$ | $\varepsilon((j,2),(i,1))$ |

Table 4-i: The 5 relevant terms of $\det \Delta_2(\varepsilon)$.

This leads to the pseudocode implementation of $SignDet\Delta_D$ shown below. It assumes that $i_0 < \ldots < i_D$ and that the sequence of subdeterminants, sorted by increasing depth, is known. The code also requires a function $SignDet_k(\Phi)$ that calculates the sign of $\det \Phi$ for a $k$-by-$k$ matrix $\Phi$. The authors have not been able to find an alternative way to determine the sign but to compute the actual determinant. Unfortunately, computing the (exact) determinant of a matrix of integers demands the use of long integer arithmetic. More about that in Section 4.3.

**function** $SignDet\Delta_D$ $(\Delta_D)$ **returns** $+1$ or $-1$
**local** $\sigma$, $k_t$, $t$
**begin**
  $t \leftarrow -1$
  **repeat**
    $t \leftarrow t + 1$
    $k_t \leftarrow k(M_t^{\Delta_D})$
    $\sigma \leftarrow SignDet_{k_t}(M_t^{\Delta_D})$
  **until** $\sigma \neq 0$
  **return** $\sigma$
**end**

Function $SignDet\Delta_D$ "scans" through the table of relevant subdeterminants. Two lines of the pseudocode, "$k_t \leftarrow k(M_t^{\Delta_D})$" and "$\sigma \leftarrow SignDet_{k_t}(M_t^{\Delta_D})$", indicate table lookups. In *Pascal* this could be implemented as a `CASE`-statement. For $D = 2$, it would consist of 5 different cases as shown below.

```
CASE t OF
  0 : s :=  SignDet2 (Pi1,Pi2,Pj1,Pj2);
  1 : s := -Sign (Pj1);
  2 : s :=  Sign (Pj2);
  3 : s :=  Sign (Pi1);
  4 : s :=  1;
END;
```

If the depth counter is of no interest, one can even unwind the loop and come up with the following code.

```
FUNCTION SignDetDelta2 (Pi1, Pi2, Pj1, Pj2): Integer;
BEGIN
  SignDetDelta2 :=  SignDet2 (Pi1, Pi2, Pj1, Pj2);
  IF SignDetDelta2 <> 0 THEN goto 999;
  SignDetDelta2 := -Sign (Pj1);
  IF SignDetDelta2 <> 0 THEN goto 999;
  :
  SignDetDelta2 :=  1;
  999: (* exit *)
END;
```

To give more insight into the computation of the terms of $\det \Delta_D(\varepsilon)$ in the order of decreasing significance, we now consider the three-dimensional case, that is,

$$\det \Delta_3(\varepsilon) = \det \begin{pmatrix} \pi_{i,1} + \varepsilon(i,1) & \pi_{i,2} + \varepsilon(i,2) & \pi_{i,3} + \varepsilon(i,3) \\ \pi_{j,1} + \varepsilon(j,1) & \pi_{j,2} + \varepsilon(j,2) & \pi_{j,3} + \varepsilon(j,3) \\ \pi_{k,1} + \varepsilon(k,1) & \pi_{k,2} + \varepsilon(k,2) & \pi_{k,3} + \varepsilon(k,3) \end{pmatrix}.$$

This polynomial has a total of 34 terms. However, only 15 of them are relevant, and those are listened in Table 2. There are two reasons why we only need to test 15 coefficients out of a total of 34. One is that the coefficient of $\varepsilon((k,3),(j,2),(i,1))$ is equal to $+1$ which is nonzero; we can therefore stop there and consider no further terms. The other reason is that certain coefficients occur more than once, that is, with different $\varepsilon$-products. For example,

$$\det \Delta_3(\varepsilon) = \cdots + \pi_{k,3} \cdot \varepsilon((j,2),(i,1)) \cdots - \pi_{k,3} \cdot \varepsilon((j,1),(i,2)) \cdots \tag{4-b}$$

Clearly, there is no need to test $-\pi_{k,3} \neq 0$, since at this depth, $+\pi_{k,3} = 0$ is already known; otherwise, the sign determination would have stopped immediately after testing the coefficient of $\varepsilon((j,2),(i,1))$.

## 4.2 Generating the Sequence of Significant Coefficients

The properly sorted sequences of $\varepsilon$-terms of the polynomials $\det \Delta_2(\varepsilon)$ and $\det \Delta_3(\varepsilon)$ are apparently very regular. In the following, this regularity will be worked out and exploited by an algorithm that automatically generates the correct sequence of $\varepsilon$-terms. This procedure can be embedded in an implementation of the function *SignDet*$\Delta_D$ that computes the sign of $\det \Delta_D(\varepsilon)$. We agree that a procedure that generates each term of $\det \Delta_D(\varepsilon)$ by collecting the proper rows and columns of the original matrix is, in a practical sense, much slower than a straight-line program that scans through a fixed sequence of submatrices. However, in higher dimensions the former might be the better strategy, since the likelihood of $\det M_\tau^{\Delta_D} = 0$ for all $\tau$ with $0 \leq \tau \leq t$ decreases very fast as $t$ increases. Let alone the fact that the tables of relevant terms for $\det \Delta_D(\varepsilon)$ becomes rather long for large $D$. The algorithm to be described can also be used for automatic generation of such tables and even for the automatic generation of codes implementing them.

We now discuss in detail how we can extract the individual terms of the polynomial $\det \Delta_D(\varepsilon)$. Recall that a term is of the form $b \cdot \varepsilon^c$, where $b$ is called the coefficient and $\varepsilon^c$ is the $\varepsilon$-product of the term. If $\varepsilon^c = \varepsilon((i_1,j_1),\ldots,(i_k,j_k))$ (so it is a $k$-fold $\varepsilon$-product), then we call $\varepsilon(i_\iota,j_\iota)$ *active*, for $1 \leq \iota \leq k$. Given the $\varepsilon$-product of a term we can extract the coefficient $b$ from the given matrix by crossing out all rows and columns that contain an active $\varepsilon(i_\iota, j_\iota)$. In order to avoid extensive double indexing and index inversions, we assume that the points whose coordinates are the entries in the $D$ rows of the matrix $\Delta_D$ have indices 1 through $D$. This allows us to ignore the difference between a point index and the corresponding row index. Indeed, this assumption is no loss of generality since the only property used in computing the sign of $\det \Delta_D(\varepsilon)$ is that the point indices are sorted and therefore the actual values are irrelevant. With this assumption, $\varepsilon(i_\iota, j_\iota)$ is in the $i_\iota$-th row and the $j_\iota$-th column and we cross out rows $i_1, i_2, \ldots, i_k$ and columns $j_1, j_2, \ldots, j_k$. This leaves a $(D-k)$-by-$(D-k)$ submatrix. Table 2 illustrates these definitions for $D = 3$. If $b \cdot \varepsilon^c$ is the term of depth $t$, then the notation in Table 2 is such that $b = \det M_t^{\Delta_D}$, $\varepsilon^c = \varepsilon_t$, and $k_t$ is the number of rows (or columns) of $M_t^{\Delta_D}$.

| $t$ | $k_t \cdot k_t$ | $v_t$ | $\det M_t^{\Delta_3}$ | $\varepsilon_t$ |
|---|---|---|---|---|
| 0 | $3 \cdot 3$ | $[4,4,4;4]$ | $+\det \begin{pmatrix} \pi_{i,1} & \pi_{i,2} & \pi_{i,3} \\ \pi_{j,1} & \pi_{j,2} & \pi_{j,3} \\ \pi_{k,1} & \pi_{k,2} & \pi_{k,3} \end{pmatrix}$ | $\varepsilon()$ |
| 1 | $2 \cdot 2$ | $[3,4,4;4]$ | $+\det \begin{pmatrix} \pi_{j,1} & \pi_{j,2} \\ \pi_{k,1} & \pi_{k,2} \end{pmatrix}$ | $\varepsilon(i,3)$ |
| 2 | $2 \cdot 2$ | $[2,4,4;4]$ | $-\det \begin{pmatrix} \pi_{j,1} & \pi_{j,3} \\ \pi_{k,1} & \pi_{k,3} \end{pmatrix}$ | $\varepsilon(i,2)$ |
| 3 | $2 \cdot 2$ | $[1,4,4;4]$ | $+\det \begin{pmatrix} \pi_{j,2} & \pi_{j,3} \\ \pi_{k,2} & \pi_{k,3} \end{pmatrix}$ | $\varepsilon(i,1)$ |
| 4 | $2 \cdot 2$ | $[3,3,4;4]$ | $-\det \begin{pmatrix} \pi_{i,1} & \pi_{i,2} \\ \pi_{k,1} & \pi_{k,2} \end{pmatrix}$ | $\varepsilon(j,3)$ |
| 5 | $1 \cdot 1$ | $[2,3,4;4]$ | $+\det(\pi_{k,1}) = +\pi_{k,1}$ | $\varepsilon((j,3),(i,2))$ |
| 6 | $1 \cdot 1$ | $[1,3,4;4]$ | $-\det(\pi_{k,2}) = -\pi_{k,2}$ | $\varepsilon((j,3),(i,1))$ |
| 7 | $2 \cdot 2$ | $[2,2,4;4]$ | $+\det \begin{pmatrix} \pi_{i,1} & \pi_{i,3} \\ \pi_{k,1} & \pi_{k,3} \end{pmatrix}$ | $\varepsilon(j,2)$ |
| 8 | $1 \cdot 1$ | $[1,2,4;4]$ | $+\det(\pi_{k,3}) = +\pi_{k,3}$ | $\varepsilon((j,2),(i,1))$ |
| 9 | $2 \cdot 2$ | $[1,1,4;4]$ | $-\det \begin{pmatrix} \pi_{i,2} & \pi_{i,3} \\ \pi_{k,2} & \pi_{k,3} \end{pmatrix}$ | $\varepsilon(j,1)$ |
| 10 | $2 \cdot 2$ | $[3,3,3;4]$ | $+\det \begin{pmatrix} \pi_{i,1} & \pi_{i,2} \\ \pi_{j,1} & \pi_{j,2} \end{pmatrix}$ | $\varepsilon(k,3)$ |
| 11 | $1 \cdot 1$ | $[2,3,3;4]$ | $-\det(\pi_{j,1}) = -\pi_{j,1}$ | $\varepsilon((k,3),(i,2))$ |
| 12 | $1 \cdot 1$ | $[1,3,3;4]$ | $+\det(\pi_{j,2}) = +\pi_{j,2}$ | $\varepsilon((k,3),(i,1))$ |
| 13 | $1 \cdot 1$ | $[2,2,3;4]$ | $+\det(\pi_{i,1}) = +\pi_{i,1}$ | $\varepsilon((k,3),(j,2))$ |
| 14 | $0 \cdot 0$ | $[1,2,3;4]$ | $+\det() = +1$ | $\varepsilon((k,3),(j,2),(i,1))$ |

Table 4-ii: The 15 relevant terms of $\det \Delta_3(\varepsilon)$.

Note that we did not yet specify how we can decide whether $b$ is $-1$ or $+1$ times the determinant of the submatrix. We now describe a rule that is based on the number of transpositions needed to sort a certain permutation. For row $\iota$, $1 \le \iota \le D$, let $j_\iota$ be the column such that $\varepsilon(\iota, j_\iota)$ is active in the term that we currently consider. By definition of a determinant there can be at most one such column but it could very well be that there is no such column. In this case we choose $j_\iota$ such that $\pi_{\iota,j_\iota}$ belongs to the main diagonal of the submatrix that was obtained after crossing out rows and columns as described above. If the number of exchanges needed to sort $(j_1, j_2, \ldots, j_D)$ is odd, then $b = \det M_t^{\Delta_D}$ is $-1$ times the determinant of the submatrix, otherwise, it is $+1$ times this determinant.

Interestingly, the number of exchanges needed to sort the sequence $(j_1, j_2, \ldots, j_D)$ is even if and only if $i_\iota + j_{i_\iota}$ is odd for an even number of pairs $(i_\iota, j_{i_\iota})$, $1 \le \iota \le k$. To see this notice that the total number of pairs $(\kappa, j_\kappa)$ with $\kappa + j_\kappa$ odd is even since

$$\sum_{\kappa=1}^{D}(\kappa + j_\kappa) = 2\sum_{\kappa=1}^{D}\kappa.$$

Now observe that $(j_1, j_2, \ldots, j_D)$ can be sorted using only exchanges of adjacent columns, that is, of integers $j_\kappa$ that differ by one. Note also that we can dispense with all exchanges between two columns where both contain an active $\varepsilon(i,j)$ or both do not. Thus, every exchange of two columns increases or decreases the number of pairs $(i_\iota, j_{i_\iota})$ with $i_\iota + j_{i_\iota}$ odd by one, which implies the claim. This property will be used in the algorithm that computes the proper sign.

The key observation that allows us to automatically generate the relevant terms of $\det \Delta_D(\varepsilon)$ is that $\varepsilon((i_1,j_1),\ldots,(i_k,j_k))$ is the $\varepsilon$-product of a relevant term if and only if $i_1 < \ldots < i_k$ and $j_1 < \ldots < j_k$. In other words, the $\varepsilon(i_\iota, j_\iota)$ go monotonically from the left top to the right bottom of the matrix. To see this take an $\varepsilon$-product that does not satisfy this condition and consider the $\varepsilon$-product defined by the same $2k$ indices that is obtained by matching the smallest $i_\iota$ with the smallest $j_\iota$, the two second smallest indices, etc. This new $\varepsilon$-product is more significant than the old one since the exponent of $\varepsilon$ it defines is smaller than the exponent of the old $\varepsilon$-product. Furthermore, the coefficients that correspond to the two $\varepsilon$-products have the same absolute value, namely the determinant of the submatrix obtained by crossing out rows $i_\iota$ and columns $j_\iota$, for $1 \le \iota \le k$.

The algorithm that generates the $\varepsilon$-products and their corresponding coefficients uses a vector

$$v = [v_1, \ldots, v_D; v_{D+1}]$$

where each $v_i$ is an integer between 1 and $D+1$ and $v_i$ corresponds to the $i$-th row of $\det \Delta_D(\varepsilon)$; $v_{D+1}$ is set equal to $D+1$ and is used only for convenience. The interpretation of $v$ is as follows: To encode the $\varepsilon$-product $\varepsilon((i_1,j_1),\ldots,(i_k,j_k))$ we set $v_{i_\iota} = j_\iota$ for $1 \le \iota \le k$. For every $i$ such that the $i$-th row does not contain an active $\varepsilon(i,j)$, we define $v_i = v_{i_\iota}$ with $i_\iota$ the smallest integer in $\{i_1, \ldots, i_k, D+1\}$ that is larger than $i$. Thus, $v_\kappa$ in $v$ implies that $\varepsilon(\kappa, v_\kappa)$ is active if and only if $v_\kappa < v_{\kappa+1}$. For example $v = [3,4,4;4]$ implies that the $\varepsilon$-product of the encoded term is $\varepsilon(1,3)$. Other examples can be found in Table 2, which gives the vectors of all relevant terms in $\det \Delta_3(\varepsilon)$.

The next problem that we face is how to generate the terms of $\det \Delta_D$ in the correct order, that is, in the order of decreasing significance. Here we use the fact that $v = [v_1, \ldots, v_d; v_{D+1}]$ encodes a more significant term than $v' = [v'_1, \ldots, v'_d; v'_{D+1}]$ if and only if $v_j > v'_j$ for $j$ the largest index, such that $v_j \ne v'_j$. This implies that $v = [D+1, \ldots, D+1; D+1]$ encodes the most significant term and, indeed, it encodes $\varepsilon() = 1$, whose coefficient is the determinant of the entire original matrix. It is now easy to write a function that computes for a given vector its successor.

```
function Next_v (v) returns Vector
local ι, κ
begin
   ι ← 1
   while v_ι = 1 do ι ← ι + 1
   v_ι ← v_ι − 1
   for κ ← ι − 1 down to 1 do v_κ ← v_ι
   return v
end
```

The alert reader will have noticed that the above function returns an "illegal" vector if the input vector is $[1, \ldots, 1; D+1]$ which is not a problem because the determinant evaluation is such that already $[1, 2, \ldots, D; D+1]$ encodes a nonzero coefficient thus there is no reason to call *Next_v* again.

After initializing $v$ to $[D+1, \ldots, D+1; D+1]$, successive calls to *Next_v* give the desired sequence of vectors. It remains to be shown how the coefficient of the encoded term can be computed. The procedure below decodes $v$ and returns the submatrix $M$ obtained after deleting the proper rows and columns from $\Delta_D$. It also returns $s$ equal to $-1$ or $+1$ depending on whether the coefficient equals $-\det M$ or $+\det M$, and returns $k$ which is equal to the number of rows (or columns for that matter) of $M$.

```
procedure Matrix (v, s, k, M)
global Δ_D, D
local ι
begin
   M ← Δ_D
   k ← D
   s ← +1
   for ι ← 1 to d do
      if v_ι < v_{ι+1} then
         {in this case ε(ι, v_ι) is active}
         if odd(ι + v_ι) then s ← −s
         delete row ι from M
         delete column v_ι from M
         k ← k − 1
end
```

We can now modify the code of *SignDet*$\Delta_D$ by replacing the table lookup by appropriate calls to *Next_v* and *Matrix*. With additional modifications the same algorithm can be used to generate the table of relevant terms in $\det \Delta_D(\varepsilon)$ or even to generate the corresponding code for *SignDet*$\Delta_D$ for any $D$. Note that, in the latter case, the loop in *SignDet*$\Delta_D$ is to be repeated only until $k_t = 0$, since in "generating mode" the values of the determinants are not computed and thus there is no natural abortion of the cycle of calls. The result for $D = 4$ can be seen in Table 6 in the Appendix.

A nice feature of the above algorithms is that we only need to change the initialization of $v$ to $[D, \ldots, D; D]$ to get an implementation for *SignDet*$\Lambda_D$ which computes the sign of the $\varepsilon$-polynomial $\det \Lambda(\varepsilon)$. For this case, the loop over all relevant terms has to be repeated until either the corresponding cofactor is nonzero, or, if we are in "generating mode," until $k_t = 1$. See Tables 3, 4 and 5 in the Appendix for the relevant terms of $\det \Lambda_D(\varepsilon)$ for $D = 2, 3, 4$. It seems worthwhile to mention that Cartesian coordinates should be used whenever possible. This reduces the problem roughly by one "dimension" compared to the homogeneous case (compare e.g., Tables 2 and 4).

The presented $\varepsilon$-polynomials $\det \Delta_D(\varepsilon)$ and $\det \Lambda_D(\varepsilon)$ illustrate that the computational overhead caused by SoS is acceptably small. One has to keep in mind that the most significant term of these $\varepsilon$-determinants corresponds to the original determinant which expresses the primitive. So, there is no way around the evaluation of the sign of this determinant for any implementation. If the input data is nondegenerated the cost of SoS is obviously zero and, in general, it is rather unlikely that the polynomials have to be evaluated down to large depths. Indeed, the largest depth or the sum of all depths that occurs in a computation can be used as a measure for the degree of degeneracy of the input data.

By evaluating the subdeterminants we systematically take care of all possible degenerate cases. Take for example the evaluation of $\det \Lambda_3(\varepsilon)$. Different cases can be distinguished by looking at the largest depth $t_{\max}$ reached during the computations. This $t_{\max}$ can be 0, 1, 2, 3, or 4 and the corresponding degeneracy is as follows (compare with Table 4 in the Appendix):

- $t_{\max} = 0$: The three points $p_i$, $p_j$, and $p_k$ are in general position.
- $t_{\max} = 1$: The three points are collinear but $p_j \neq p_k$ and the line containing the three points is not vertical.
- $t_{\max} = 2$: The three points lie on a common vertical line but $p_j \neq p_k$.
- $t_{\max} = 3$: Point $p_j$ coincides with $p_k$, but not with $p_i$, and the line through $p_i$ and $p_j$ is not vertical.
- $t_{\max} = 4$: All three points lie on a common vertical line and, $p_j = p_k$.

It would be interesting to see this somewhat unnatural case analysis in greater detail since it gives a nonobvious breakdown into degenerate cases that has curious properties.

This discussion completes the implementation of SoS with respect to the predicate $Positive_d$ for point sets in $E^d$. We considered both the Cartesian and the homogeneous case. The key was to find a method that generates the proper sequence of relevant terms of $\det \Delta_D(\varepsilon)$ and $\det \Lambda_D(\varepsilon)$ ordered by decreasing significance. With this, the implementation of the functions $SignDet\Delta_D$ and $SignDet\Lambda_D$ was easy. We will see in Section 5 that both functions can also be used to implement other predicates.

## 4.3 Remarks on the Sign Computation of Determinants

In the previous sections we reduced all computations to a sequence of sign evaluations of determinants. In the primitives discussed in this paper, the matrices are at most of size $(d+2)$-by-$(d+2)$, $d$ the dimension of the space, and all elements are assumed to be integers. Theoretically, the sign of such a determinant can be determined in constant time if we assume that $d$ is a constant. This assumption is indeed fair since SoS is intended primarily for low-dimensional geometric computations. In practice, however, it is important to optimize the sign computation since it will be in the innermost loop of every program that uses SoS — which does not mean that this issue is less important for programs not employing SoS. We remark on a few methods that can be used to get speed in these computations.

One important condition that we have to meet is that the sign of the determinant has to be computed exactly — we cannot tolerate a +1 for a 0, etc. Assuming that the coordinates or parameters are integers, we can either use long integer arithmetic or modular arithmetic based on the Chinese remainder theorem. For details on both methods refer to [Kn69]. If we actually compute the determinant in order to find its sign — and no method is known to the authors that avoids the actual computation of the determinant — we have to be prepared to deal with numbers of absolute size at least $\mu^D$, where $\mu$ denotes maximum absolute value of any data item and $D$ denotes the largest size of matrices we work with. To see this, just take the $D$-by-$D$ matrix whose entries are all zero except for the ones in the main diagonal, where they are equal to $\mu$; the determinant of this matrix is $\mu^D$. An upper bound on the absolute value of the determinants is given by a well known theorem of Hadamard that states that

$$|\det(\pi_{1\ldots D,1\ldots D})| \leq \prod_{i=1}^{D} \sqrt{\sum_{j=1}^{D} \pi_{i,j}^2} \leq \mu^D D^{\frac{D}{2}}.$$

Among other things this upper bound on the absolute value of a determinant gives us an upper bound on the number of computer words needed for the computation if we use long integer arithmetic.

Without any hardware support long integer arithmetic is very time consuming, which might motivate us to resort to the use of approximation methods. Any computation of the determinant using floating-point arithmetic of bounded length is such an approximation. Floating-point arithmetic is usually rather fast since it enjoys the needed hardware support on most of today's computers. If the value that we get is sufficiently far from zero, we can be sure that the correct value is different from zero and lies on the same side of zero. But how can we quantify "sufficiently far from zero"? In any case, we could now use Gaussian elimination (see e.g., [GVL83]) which takes $O(D^3)$ time or asymptotically faster methods based on matrix multiplication as described for instance in [AHU74]. We do not believe that the latter methods could be of any practical use, though. However, if the value that we get is suspiciously close to zero, we have to use some other method to determine the sign of the determinant.

Finally, we would like to mention that the determinant of a $D$-by-$D$ matrix can be expressed in terms of subdeterminants, and that some of these subdeterminants might later appear again when the evaluation of $\det \Delta_D(\varepsilon)$ or $\det \Lambda_D(\varepsilon)$ proceeds. It is conceivable that the values of such subdeterminants are saved and used again when needed. Even so, we do not believe that such a method could lead to significant savings since we expect that on the average only very few terms of the $\varepsilon$-determinants are needed.

# 5 Further Applications of SoS for Determinants

In this section we demonstrate that the algorithmic solution to many geometric problems can be based on primitive operations that compute the sign of determinants. Those include problems that deal with objects different from points. There are two major reasons why determinants are useful beyond problems for points. One is that more complicated geometric objects are often given by a finite set or sequence of points. Examples are line segments given by two points and triangles specified by three points. This will be illustrated in Section 5.1, which revisits the Parity Algorithm discussed in the Introduction. The other reason (and this is the more profound although less obvious of the two) is that other objects can be thought of as points in a different space. Take, for example, a hyperplane in $d$ dimensions. It can be specified by a linear relation of the form

$$\eta_1 x_1 + \eta_2 x_2 + \cdots + \eta_d x_d + \eta_{d+1} = 0.$$

Multiplying the above relation with a nonzero constant does not change the hyperplane. This suggests that we think of the hyperplane as the point with homogeneous coordinates

$$(\eta_1, \eta_2, \ldots, \eta_d; \eta_{d+1})$$

in $d$ dimensions. This view of hyperplanes will be discussed in more detail in Sections 5.2 and 5.3. Of course, an $n$-gon specified by a sequence of $n$ points in the plane can be interpreted as a point too — in this case it is a point in $2n$ dimensions. However, in contrast to the former case, this view is not likely to lead to any useful application of determinants since it becomes increasingly expensive to compute them as the size of the matrix increases. Finally, Section 5.4 shows that even nonlinear geometric objects such as circles and spheres can profitably be interpreted as points in low dimensions as well.

By no means do we believe that the list of applications for primitives concerning the sign of determinants, as presented in this paper, is exhaustive. In fact, because of the versatility of determinants, an enumeration of their applications in geometric computation is far beyond the scope of this paper. We agree though that such an enumeration is a challenging task.

## 5.1 Point-in-Polygon Test

Recall the Parity Algorithm for the point-in-polygon problem sketched in the Introduction. In order to test whether a given point $p$ lies inside a simple polygon $P$, the algorithm intersects the horizontal half-line $r$, whose left endpoint is $p$, with all edges of polygon $P$. If the number of edges intersecting $r$ is odd, then $p$ lies inside $P$, and if this number is even, $p$ lies outside. The subtlety of this algorithm lies in the treatment of special cases since the above characterization holds, in general, only if we introduce certain artificial counting mechanisms whenever $r$ contains a vertex or even an entire edge of $P$. In this section we show that the test whether or not an edge intersects the horizontal half-line $r$ can be reduced to computing the signs of certain determinants. SoS is then used to simulate a perturbation of the point and the polygon which removes all degeneracies. The algorithm assumes that $P$ is given by a sequences of vertices $(v_1, v_2, \ldots, v_n)$ and that all coordinates including those of $p$ are integers.

We consider now the problem to test whether $r$ intersects an edge $e$ of $P$ given by its two endpoints. Let $u = (v_1, v_2)$ and $w = (\omega_1, \omega_2)$ be the two endpoints and recall that $p = (\pi_1, \pi_2)$ is the left endpoint of $r$. Because of SoS we can assume that $u$, $w$, $p$ are not collinear and that no two of the three points lie on a common horizontal line. Note first that $r$ and $e$ intersect only if the second coordinate of $p$ lies between the second coordinates of $u$ and $w$. Assume $v_2 < \omega_2$. If indeed $v_2 < \pi_2 < \omega_2$ then $r \cap e \neq \emptyset$ if and only if $(u, w, p)$ defines a left-turn (see Figure 5-I).

It is now not very difficult to develop this case analysis into a predicate that tests for intersection. To perturb the points we use the same $\varepsilon$-expansion as described in Section 3.2, that is, we replace $v_i = (\nu_{i,1}, \nu_{i,2})$ by $v_i(\varepsilon) = (\nu_{i,1}(\varepsilon), \nu_{i,2}(\varepsilon))$ where $\nu_{i,j}(\varepsilon) = \nu_{i,j} + \varepsilon(i,j)$ with $\varepsilon(i,j)$ as in (3-b). For a uniform treatment we define $p = v_0 = (\nu_{0,1}, \nu_{0,2})$ and write the predicate for arbitrary three vertices rather than for $v_0$ and two successive vertices of $P$.

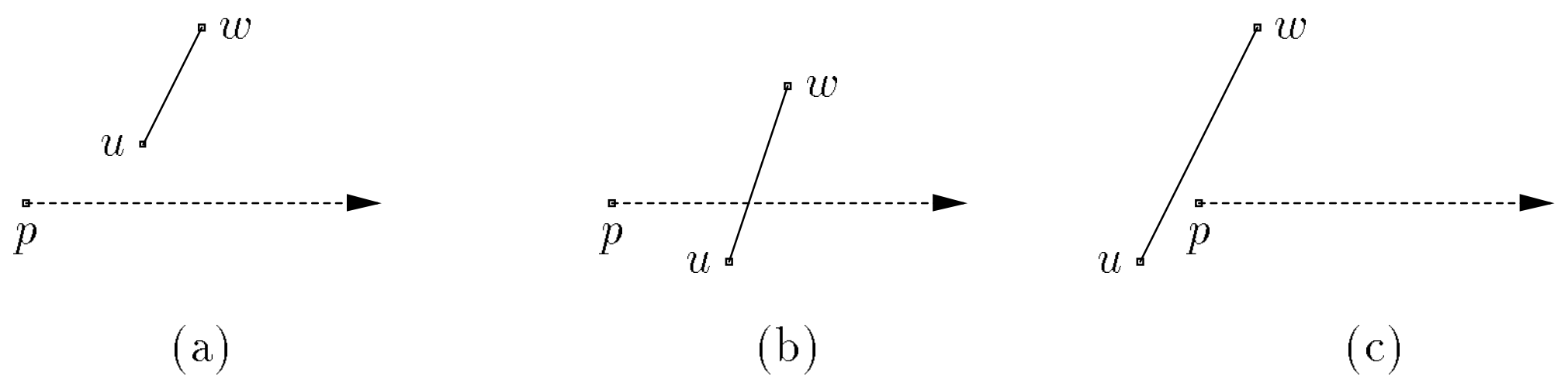

Figure 5-I: The three cases to consider for $r \cap e$ using SoS.

In (a), $r$ and $e$ do not intersect since the second coordinate of $p$ does not lie between those of $u$ and $w$. In (c), they do not intersect since $(u, w, p)$ is a right-turn.

**Predicate 3** (*IntersectHalfLine*) Let $v_i$, $v_j$, and $v_k$ be three vertices with pairwise different indices $0 \leq i, j, k \leq n$. The following pseudocode returns true if the edge from $v_j(\varepsilon)$ to $v_k(\varepsilon)$ intersects the horizontal half-line whose left endpoint is given by $v_i(\varepsilon)$, and false otherwise.

**function** *IntersectHalfLine* $(v_i; v_j, v_k)$ **returns** Boolean
**local** $i'$, $j'$, $k'$, $s'$, $d'$
**begin**
  W.l.o.g. assume $Smaller(\nu_{j,2}; \nu_{k,2})$.
  **if** $Smaller(\nu_{j,2}; \nu_{i,2}) \wedge Smaller(\nu_{i,2}; \nu_{k,2})$ **then**
    $Sort_3((i,j,k), (i',j',k'), s')$

$$d' \leftarrow SignDet\Lambda_3 \begin{pmatrix} \nu_{i',1}(\varepsilon) & \nu_{i',2}(\varepsilon) & 1 \\ \nu_{j',1}(\varepsilon) & \nu_{j',2}(\varepsilon) & 1 \\ \nu_{k',1}(\varepsilon) & \nu_{k',2}(\varepsilon) & 1 \end{pmatrix}$$

    **if** odd$(s')$ **then** $d' \leftarrow -d'$
    **return** $(d' = +1)$
  **else**
    **return** false
**end**

A few remarks are in order. When the above function is applied to the point-in-polygon problem, $i = 0$ always holds. Thus, the sorting of $(i, j, k)$ can be reduced to a single comparison between $j$ and $k$. Furthermore, to avoid all degeneracies for the point-in-polygon test, it is sufficient to perturb only the point $p = v_0$. Indeed, if

$$\det \begin{pmatrix} \nu_{i,1}(\varepsilon) & \nu_{i,2}(\varepsilon) & 1 \\ \nu_{j,1} & \nu_{j,2} & 1 \\ \nu_{k,1} & \nu_{k,2} & 1 \end{pmatrix} = 0,$$

then we necessarily have $\nu_{j,2} = \nu_{k,2} \neq \nu_{i,2}(\varepsilon)$ and therefore, the determinant does not even get evaluated. The savings that one gets this way are only nominal which we interpret as an argument for the efficiency of our general method.

The remainder of this section is used to comment on what happens if the test point $p$ lies on the boundary of the polygon $P$. If we use the above primitive as is, SoS will neglect this special case and find that $p$ lies on either side of $P$'s boundary. The decision depends on the relative positions of $p$ and the vertices of $P$, and we might as well assume that it is arbitrary although consistent.

Such a decision may or may not be desirable. If it is not acceptable, one could test whether or not $p$ lies on the boundary of $P$ before running the Parity Algorithm with SoS. Once more this test can be reduced to computing signs of determinants.

## 5.2 Hyperplanes in Euclidean Space

Algorithms for hyperplanes play a central role in computational geometry. This becomes obvious when one thinks of the importance of problems such as linear programming, computing the intersection of half-spaces, and constructing arrangements of hyperplanes (see [PS85] and [Ed87] for further details and references). The goal of this section is to demonstrate how the techniques of Section 4 can be used to implement a typical primitive operation needed in those algorithms. This will open up an entire class of problems to the use of SoS. The main tool that lets us exploit the techniques of Section 4 when we handle hyperplanes is a duality transformation that maps hyperplanes to points and vice versa. In essence, this transform is nothing but a reinterpretation of what hyperplanes and points are.

In this section we assume that a hyperplane $h$ in $E^d$ is specified by its nonzero normal vector $a = (\alpha_1, \ldots, \alpha_d)$ and a number, $-\alpha_{d+1}$, called the *offset*. Now, a hyperplane $h$ consists of all points $x \in E^d$ such that

$$\langle x, a \rangle + \alpha_{d+1} = 0, \qquad \text{(5-a)}$$

that is, the scalar product of $x$ and $a$ equals the offset. Notice that the hyperplane does not change if we multiply the normal vector and the offset by some nonzero number. We define $h^*$ as the point whose homogeneous coordinates are $(\alpha_1, \ldots, \alpha_d; \alpha_{d+1})$. Geometrically speaking, $h^*$ lies on the line through the origin defined by $a$, and the distance of $h^*$ from the origin is the inverse of the distance between $h$ and the origin. This can easily verified after observing that $|\alpha_{d+1}|$ is the distance between $h$ and the origin, provided $a$ has unit length. Note also that the origin lies between $h$ and $h^*$ (see Figure 5-II). Conversely, for a point $p$ with homogeneous coordinates $(\pi_1, \pi_2, \ldots, \pi_d; \pi_{d+1})$ we let $p^*$ be the hyperplane with normal vector $(\pi_1, \pi_2, \ldots, \pi_d)$ and offset $-\pi_{d+1}$.

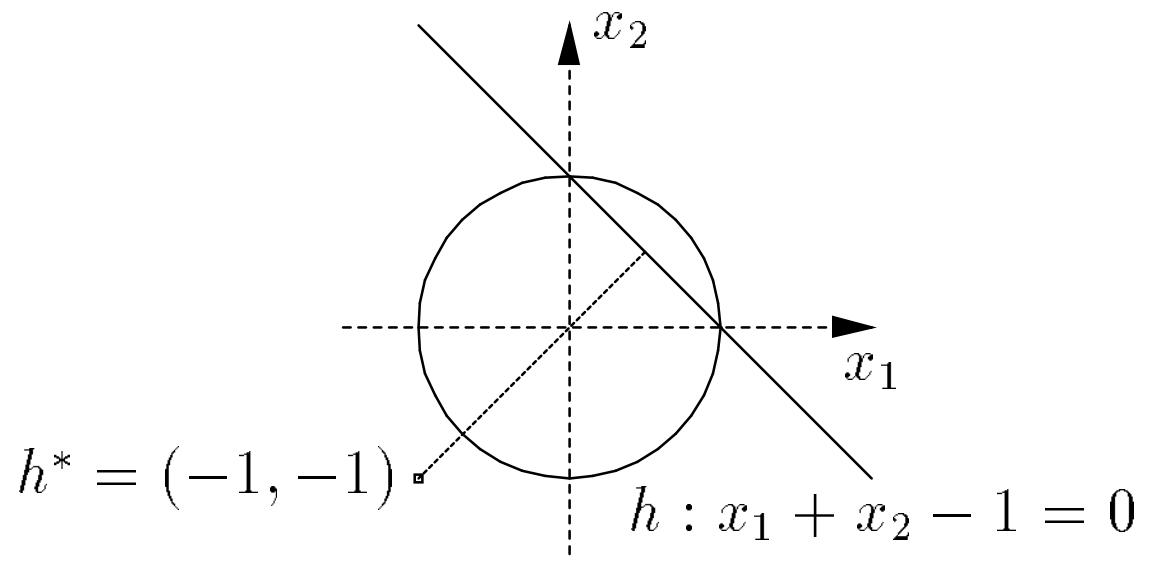

Figure 5-II: Mapping a line to a point and vice versa.

It is straightforward to show that this transformation preserves incidences; that is, $p \in h$ if and only if $h^* \in p^*$. Indeed, it is a triviality when one remembers what $p \in h$ means algebraically, namely, that

$$\pi_1\alpha_1 + \pi_2\alpha_2 + \cdots + \pi_d\alpha_d + \pi_{d+1}\alpha_{d+1} = 0.$$

It is equally easy to prove that this mapping preserves the relative order between a point and a hyperplane. To describe what exactly we mean by this define

$$h^+ = \{x | \langle x, a \rangle + \alpha_{d+1} > 0\} \quad \text{and} \quad h^- = \{x | \langle x, a \rangle + \alpha_{d+1} < 0\},$$

and call those the *positive* and *negative sides* or *half-spaces* of $h$. By order preservation we mean that $p \in h^+$ if and only if $h^* \in p^{*+}$. Here, a warning is appropriate to avoid future confusion. If we multiply the normal vector and the offset of a hyperplane $h$ by $-1$ we do not change the hyperplane but we do change the sides of $h$; what was previously its positive side is now its negative one and the other way round. We will take advantage of this curiosity by encoding the positive and negative sides into the hyperplane's specification. Note that geometrically, the normal vector of a hyperplane $h$ points to its positive side.

The primitive operation that we wish to tackle in this section is to decide on which side of a hyperplane $h_{i_d}$ the intersection of $d$ other hyperplanes $h_{i_0}, \ldots, h_{i_{d-1}}$ lies. By the use of SoS, the absence of any kind of degeneracies can be assumed; so $h_{i_0}(\varepsilon)$ through $h_{i_{d-1}}(\varepsilon)$ intersect in a unique point which does not lie on $h_{i_d}(\varepsilon)$. By Cramer's rule, the intersection point $p = (\pi_1, \pi_2, \ldots, \pi_d)$ of $d$ hyperplanes is given by the coordinates

$$\pi_i = \frac{\det \Delta_{d,i}}{\det \Delta_d},$$

where $\Delta_d$ is the matrix

$$\begin{pmatrix} \alpha_{i_0,1} & \alpha_{i_0,2} & \cdots & \alpha_{i_0,d} \\ \alpha_{i_1,1} & \alpha_{i_1,2} & \cdots & \alpha_{i_1,d} \\ \vdots & \vdots & \ddots & \vdots \\ \alpha_{i_{d-1},1} & \alpha_{i_{d-1},2} & \cdots & \alpha_{i_{d-1},d} \end{pmatrix},$$

and $\Delta_{d,i}$ is the same matrix after replacing the $i$-th column from the left by the vector

$$\begin{pmatrix} -\alpha_{i_0,d+1} \\ -\alpha_{i_1,d+1} \\ \vdots \\ -\alpha_{i_{d-1},d+1} \end{pmatrix}.$$

Point $p$ lies in the positive half-space of $h_d$ if and only if

$$\pi_1\alpha_{i_d,1} + \pi_2\alpha_{i_d,2} + \cdots + \pi_d\alpha_{i_d,d} + \alpha_{i_d,d+1} > 0.$$

Provided that $\det \Delta_d$ is positive, this is equivalent to

$$\det \Delta_{d,1}\alpha_{i_d,1} + \det \Delta_{d,2}\alpha_{i_d,2} + \cdots + \det \Delta_{d,d}\alpha_{i_d,d} + \det \Delta_d\alpha_{i_d,d+1} > 0.$$

In case of a negative $\det \Delta_d$ the above statement is valid after reversing the direction of the inequality. Consequently, $p \in h_d^+$ if and only if

$$\det \Delta_{d+1} \cdot \det \Delta_d > 0.$$

This can be seen by developing

$$\det \Delta_{d+1} = \det \begin{pmatrix} \alpha_{i_0,1} & \alpha_{i_0,2} & \cdots & \alpha_{i_0,d} & \alpha_{i_0,d+1} \\ \alpha_{i_1,1} & \alpha_{i_1,2} & \cdots & \alpha_{i_1,d} & \alpha_{i_1,d+1} \\ \vdots & \vdots & \ddots & \vdots & \vdots \\ \alpha_{i_{d-1},1} & \alpha_{i_{d-1},2} & \cdots & \alpha_{i_{d-1},d} & \alpha_{i_{d-1},d+1} \\ \alpha_{i_d,1} & \alpha_{i_d,2} & \cdots & \alpha_{i_d,d} & \alpha_{i_d,d+1} \end{pmatrix}$$

using the last row. Now, we can use this to write a procedure that decides on which side of a hyperplane $d$ other hyperplanes intersect. It uses SoS, as described in Section 3.3.

**Predicate 4** (*OnPositiveSide*) Let $h_{i_0}, h_{i_1}, \ldots, h_{i_d}$ be $d+1$ hyperplanes in $d$ dimensions, given as in (5-a), and with distinct indices $0 \leq i_0, i_1, \ldots, i_d \leq n-1$. The following function written in pseudocode returns true if the intersection $\bigcap_{\lambda=0}^{d-1} h_{i_\lambda}(\varepsilon)$ lies in the positive half-space of $h_{i_d}(\varepsilon)$ and false if it lies in the negative one.

**function** $OnPositiveSide_d(h_{i_0}, \ldots, h_{i_{d-1}}; h_{i_d})$ **returns** Boolean
**local** $i'_0, \ldots, i'_{d-1}, s', d', i''_0, \ldots, i''_{d-1}, i''_d, s'', d''$
**begin**
$Sort_d((i_0, \ldots, i_{d-1}), (i'_0, \ldots, i'_{d-1}), s')$
$Sort_{d+1}((i_0, \ldots, i_{d-1}, i_d), (i''_0, \ldots, i''_{d-1}, i''_d), s'')$

$$d' = SignDet\Delta_d \begin{pmatrix} \alpha_{i'_0,1}(\varepsilon) & \alpha_{i'_0,2}(\varepsilon) & \cdots & \alpha_{i'_0,d}(\varepsilon) \\ \alpha_{i'_1,1}(\varepsilon) & \alpha_{i'_1,2}(\varepsilon) & \cdots & \alpha_{i'_1,d}(\varepsilon) \\ \vdots & \vdots & \ddots & \vdots \\ \alpha_{i'_{d-1},1}(\varepsilon) & \alpha_{i'_{d-1},2}(\varepsilon) & \cdots & \alpha_{i'_{d-1},d}(\varepsilon) \end{pmatrix}$$

**if** odd $(s')$ **then** $d' \leftarrow -d'$

$$d'' \leftarrow SignDet\Delta_{d+1} \begin{pmatrix} \alpha_{i''_0,1}(\varepsilon) & \alpha_{i''_0,2}(\varepsilon) & \cdots & \alpha_{i''_0,d}(\varepsilon) & \alpha_{i''_0,d+1}(\varepsilon) \\ \alpha_{i''_1,1}(\varepsilon) & \alpha_{i''_1,2}(\varepsilon) & \cdots & \alpha_{i''_1,d}(\varepsilon) & \alpha_{i''_1,d+1}(\varepsilon) \\ \vdots & \vdots & \ddots & \vdots & \vdots \\ \alpha_{i''_{d-1},1}(\varepsilon) & \alpha_{i''_{d-1},2}(\varepsilon) & \cdots & \alpha_{i''_{d-1},d}(\varepsilon) & \alpha_{i''_{d-1},d+1}(\varepsilon) \\ \alpha_{i''_d,1}(\varepsilon) & \alpha_{i''_d,2}(\varepsilon) & \cdots & \alpha_{i''_d,d}(\varepsilon) & \alpha_{i''_d,d+1}(\varepsilon) \end{pmatrix}$$

**if** odd $(s'')$ **then** $d'' \leftarrow -d''$
**return** $(d' = d'')$
**end**

## 5.3 Nonvertical Hyperplanes

In many applications we know that all hyperplanes we have to deal with are nonvertical; that is, they intersect the $d$-th coordinate axis in a unique point. Examples are Voronoi diagrams or, more generally, power diagrams for arbitrary order and weighted Voronoi diagrams (see for instance [Ed87] and [AI86]). It is beyond the scope of this paper to describe how the data for those problems are used to generate hyperplanes — it will be enough to know that they are obtained via geometric transforms which do not create vertical hyperplanes.

A nonvertical hyperplane $h$ in $d$-dimensions can be specified by a relation of the form

$$\alpha_1 x_1 + \alpha_2 x_2 + \cdots + \alpha_{d-1} x_{d-1} + x_d + \alpha_d = 0. \tag{5-b}$$

The advantage of describing a hyperplane using this form rather than the one in Section 5.2 is that it takes only $d$ parameters rather than $d+1$. This will lead to some savings when it comes to computing signs of determinants (compare for instance det $\Delta_4$ in Table 6 and det $\Lambda_4$ in Table 5). Since every hyperplane $h$ is now nonvertical, we can uniquely define what we mean when we say that a point lies (vertically) above or below $h$. Define

$$h^+ = \{x = (x_1, \ldots, x_d) | \alpha_1 x_1 + \cdots + \alpha_{d-1} + x_d + \alpha_d > 0\},$$

and let $h^- = E^d - h - h^+$. A point $p$ is said to lie *above* $h$ if $p \in h^+$ and *below* $h$ if $p \in h^-$.

The primitive operation that we consider in this section decides whether the intersection of $d$ hyperplanes $h_{i_0}, \ldots, h_{i_{d-1}}$ lies above or below hyperplane $h_{i_d}$. The use of SoS as in Section 3.2 allows us to assume that indeed $h_{i_0}$ through $h_{i_{d-1}}$ intersect in a unique point that does not lie on $h_{i_d}$. A decision procedure based on comparing the signs of the two determinants can be derived from the procedure given in Section 5.2. We just replace all $\alpha_{i,d}(\varepsilon)$ by 1 and exchange the last two columns of the second matrix in function *OnPositiveSide*$_d$. This leads to the following predicate.

**Predicate 5** (*Above*) Let $h_{i_0}, h_{i_1}, \ldots, h_{i_d}$ be $d+1$ nonvertical hyperplanes in $E^d$ specified as in (5-b), and with pairwise different indices $0 \leq i_0, i_1, \ldots, i_d \leq n-1$. The following predicate returns true if the point of intersection $\bigcap_{\lambda=0}^{d-1} h_{i_\lambda}(\varepsilon)$ lies above $h_{i_d}(\varepsilon)$, and false if it lies below $h_{i_d}(\varepsilon)$.

**function** $\mathit{Above}_d(h_{i_0}, \ldots, h_{i_{d-1}}; h_{i_d})$ **returns** Boolean
**local** $i'_0, \ldots, i'_{d-1}, s', d', i''_0, \ldots, i''_{d-1}, i''_d, s'', d''$
**begin**
  $\mathit{Sort}_d((i_0, \ldots, i_{d-1}), (i'_0, \ldots, i'_{d-1}), s')$
  $\mathit{Sort}_{d+1}((i_0, \ldots, i_{d-1}, i_d), (i''_0, \ldots, i''_{d-1}, i''_d), s'')$

$$d' = \mathit{SignDet}\Lambda_d \begin{pmatrix} \alpha_{i'_0,1}(\varepsilon) & \cdots & \alpha_{i'_0,d-1}(\varepsilon) & 1 \\ \vdots & \ddots & \vdots & \vdots \\ \alpha_{i'_{d-1},1}(\varepsilon) & \cdots & \alpha_{i'_{d-1},d-1}(\varepsilon) & 1 \end{pmatrix}$$

  **if** odd $(s')$ **then** $d' \leftarrow -d'$

$$d'' \leftarrow \mathit{SignDet}\Lambda_{d+1} \begin{pmatrix} \alpha_{i''_0,1}(\varepsilon) & \cdots & \alpha_{i''_0,d-1}(\varepsilon) & \alpha_{i''_0,d}(\varepsilon) & 1 \\ \vdots & \ddots & \vdots & \vdots & \\ \alpha_{i''_{d-1},1}(\varepsilon) & \cdots & \alpha_{i''_{d-1},d-1}(\varepsilon) & \alpha_{i''_{d-1},d}(\varepsilon) & 1 \\ \alpha_{i''_d,1}(\varepsilon) & \cdots & \alpha_{i''_d,d-1}(\varepsilon) & \alpha_{i''_d,d}(\varepsilon) & 1 \end{pmatrix}$$

  **if** odd $(s'')$ **then** $d'' \leftarrow -d''$
  **return** $(d' \neq d'')$
**end**

## 5.4 In-sphere Test

In $d$ dimensions any $d+1$ affinely independent points (i.e., points that do not lie in a common hyperplane) define a unique sphere that goes through the $d+1$ points. For example, in two dimensions there is a unique circle through any three noncollinear points. Given $d+2$ points $p_0, p_1, \ldots, p_{d+1}$ the problem we address in this section is how we can determine whether $p_{d+1}$ lies inside or outside the sphere specified by the first $d+1$ points, assuming this sphere is unique. Such a test is useful for constructing Voronoi diagrams (as shown in [GS85] for $d=2$) and other problems where circles and spheres play a role.

An elegant solution to this problem can be given using a transform that lifts a sphere in $d$ dimensions to $d+1$ dimensions where it is represented by a hyperplane. This transformation can be traced back in the literature to [Se82] and has since been used throughout the computational geometry literature (see [GS85], [PS85], and [Ed87]). For the case of circles in the plane we explain this transformation in detail and finally phrase the predicate for general dimensions.

Let $U: \; x_3 = x_1^2 + x_2^2$ be the paraboloid of revolution whose symmetry axis is the $x_3$-axis, and let

$$c: \; (x_1 - \gamma_1)^2 + (x_2 - \gamma_2)^2 = \gamma_3^2$$

be a circle in the $x_1x_2$-plane. Note that $(\gamma_1, \gamma_2)$ is the center of the circle and $\gamma_3$ is its radius. The lifting map transforms $c$ to the plane $c^*$ in three dimensions given by the equation

$$c^* : \; x_3 = 2\gamma_1 x_1 + 2\gamma_2 x_2 - (\gamma_1^2 + \gamma_2^2 - \gamma_3^2).$$

The quick reader will already have verified that the vertical projection of $U \cap c^*$, which is an ellipse in three dimensions, onto the $x_1x_2$-plane is equal to the original circle $c$. A point $p = (\pi_1, \pi_2)$ lies inside $c$ if and only if its vertical projection onto $U$ lies below $c^*$. This insight gives us some hope that, in fact, the problem can be bent such that Predicate 5 from the previous section is applicable. Before we continue our exploration in this direction, let us understand how the original statement of the problem and the lifting map are connected. Recall that there are four original points, which we call $p_0$, $p_1$, $p_2$, and $p_3$. The first three determine the circle $c$ and therefore the plane $c^*$. Moreover, if we project them vertically onto $U$, then $c^*$ is the plane through these points on the paraboloid $U$. The question is now whether $p_3' = (\pi_{3,1}, \pi_{3,2}, \pi_{3,1}^2 + \pi_{3,2}^2)$, which is the vertical projection of $p_3$ onto $U$, lies below $c^*$ (in which case $p_3$ lies inside $c$) or above $c^*$ (then $p_3$ lies outside $c$). By the use of SoS, we can assume that the four points are in general position.

This problem can be mapped to the plane problem of the previous section if we use a dual transform. This transform replaces each point on $U$ by the unique plane whose intersection with $U$ is this point. If $p_3' = (\pi_1, \pi_2, \pi_1^2 + \pi_2^2)$, then the formula for this dual plane is

$$p^* : \; x_3 = 2\pi_1 x_1 + 2\pi_2 x_2 - (\pi_1^2 + \pi_2^2).$$

We see that this is indeed the lifting map applied to point $p = (\pi_1, \pi_2)$ in the $x_1x_2$-plane. This duality transform preserves incidences and above-below order in a way similar to the duality transform described in Section 5.2. This leaves us with the following correspondence between the original point-circle problem and the derived plane-point problem: Point $p_3$ lies inside $c$ (the circle through points $p_0$, $p_1$, and $p_2$) if and only if the intersection point of the planes $p_0^*$, $p_1^*$, and $p_2^*$ lies below $p_3^*$. The statement is also valid if we replace "inside $c$" by "outside $c$" and "below $p_3^*$" by "above $p_3^*$".

We leave the generalization of this two-dimensional exercise to three and higher dimensions to the curious reader. In any case, Predicate 5 can now be used to implement Predicate 6 which formalizes the in-sphere test in $d$ dimensions. If we apply Predicate 5 directly, we will find ourselves computing the sign of determinants of the form

$$\det \begin{pmatrix} 2\pi_{i_0,1} & \cdots & 2\pi_{i_0,d} & -(\pi_{i_0,1}^2 + \cdots + \pi_{i_0,d}^2) & -1 \\ \vdots & \ddots & \vdots & \vdots & \vdots \\ 2\pi_{i_{d+1},1} & \cdots & 2\pi_{i_{d+1},d} & -(\pi_{i_{d+1},1}^2 + \cdots + \pi_{i_{d+1},d}^2) & -1 \end{pmatrix}.$$

The sign does not change if we divide the entries in the left $d$ columns by 2. Similarly, we can remove the minus signs in the last two columns without changing the sign of the determinant. However, there remains one problem with determinants of the above type, and this is that the values in the $(d+1)$-st column from the left depend on the values in the left $d$ columns. In particular, with SoS, the $\varepsilon$-expressions of the point coordinates appear in mixed products in the $(d+1)$-st column. This turns out to be a real pain when we implement SoS for this type of determinants. A cheap trick that handles this problem is not to perturb the original points but rather to perturb the vertical projections onto the paraboloid in $d+1$ dimensions. In effect, this means that we introduce

$$\pi_{i_\lambda,d+1} = \sum_{\nu=1}^{d} \pi_{i_\lambda,\nu}^2 \quad \text{for} \quad 0 \le \lambda \le d+1, \tag{5-c}$$

and then perturb the points $(\pi_{i_\lambda,1}, \ldots, \pi_{i_\lambda,d}, \pi_{i_\lambda,d+1})$. Because the perturbation of the $(d+1)$-st coefficient does not depend on the first $d$ coefficients, this implies that the points are perturbed away from the paraboloid $U$. On the other hand, if the perturbation is small enough we are still close enough to the original situation.

**Predicate 6** (*InSphere*) Let $p_{i_0}, p_{i_1}, \ldots, p_{i_{d+1}}$ be $d+2$ points in $d$ dimensions with pairwise different indices in the range from 0 through $n-1$. The program below returns true if the perturbed image of $p_{i_{d+1}}$ lies inside the sphere through the perturbed images of the first $d+1$ points, and returns false if it lies outside.

**function** $\mathit{InSphere}_d\ (p_{i_0}, \ldots, p_{i_d}; p_{i_{d+1}})$ **returns** Boolean
**local** $i'_0, \ldots, i'_d,\ s',\ d',\ i''_0, \ldots, i''_d, i''_{d+1},\ s'',\ d''$
**begin**

$\mathit{Sort}_{d+1}((i_0, \ldots, i_d), (i'_0, \ldots, i'_d), s')$

$\mathit{Sort}_{d+2}((i_0, \ldots, i_d, i_{d+1}), (i''_0, \ldots, i''_d, i''_{d+1}), s'')$

$$d' = \mathit{SignDet}\Lambda_{d+1} \begin{pmatrix} \pi_{i'_0,1}(\varepsilon) & \cdots & \pi_{i'_0,d}(\varepsilon) & 1 \\ \vdots & \ddots & \vdots & \vdots \\ \pi_{i'_d,1}(\varepsilon) & \cdots & \pi_{i'_d,d}(\varepsilon) & 1 \end{pmatrix}$$

**if** odd $(s')$ **then** $d' \leftarrow -d'$

Set $\pi_{i_\lambda,d+1}$ as in (5-c).

$$d'' \leftarrow \mathit{SignDet}\Lambda_{d+2} \begin{pmatrix} \pi_{i''_0,1}(\varepsilon) & \cdots & \pi_{i''_0,d}(\varepsilon) & \pi_{i''_0,d+1}(\varepsilon) & 1 \\ \vdots & \ddots & \vdots & \vdots & \\ \pi_{i''_d,1}(\varepsilon) & \cdots & \pi_{i''_d,d}(\varepsilon) & \pi_{i''_d,d;1}(\varepsilon) & 1 \\ \pi_{i''_{d+1},1}(\varepsilon) & \cdots & \pi_{i''_{d+1},d}(\varepsilon) & \pi_{i''_{d+1},d+1}(\varepsilon) & 1 \end{pmatrix}$$

**if** odd $(s'')$ **then** $d'' \leftarrow -d''$

**return** $(d' = d'')$

**end**

Note that the rightmost column of the first matrix in the above program should really consist of $-1$'s. To stress the similarity with predicate $\mathit{Above}_{d+1}$ in the previous section we replaced the $-1$'s by $+1$'s and thus changed the sign of $d'$. This effect is compensated by the fact that we want to return true where function $\mathit{Above}_{d+1}$ returns false.

# 6 Remarks and Discussion

The main contribution of this paper is the introduction of a general technique that can be used to deal with degenerate input for geometric programs. The main purpose of this paper is to demonstrate that this technique (which we call SoS, the **S**imulation **o**f **S**implicity) is immanently practical, despite its high-powered appearance. Indeed, the authors believe that SoS will become a standard tool for implementing geometric algorithms. A pragmatic consequence of this technique is that authors of geometric algorithms can now be more confident about the implementability of their algorithms even in the presence of any conceivable degeneracies, provided SoS is applicable to their algorithms.

This raises the question of determining the limitations of SoS — what are the properties of an algorithm that allows us to use SoS when we implement it? One important feature of algorithms that are amenable to SoS is that their algebraic computations are of constant depth. The deeper the algebraic computation the more complicated is the polynomial (or, in general, the function) in $\varepsilon$ generated by SoS, and the less tractable is its evaluation. Another limitation of SoS is the necessity of absolute precision in the evaluation of algebraic formulas. As long as square roots can be eliminated by squaring the equation and similar techniques can be used to remove other irrational functions this is not a problem, but there are cases where it is not that easy. Typical examples for such problem cases are algorithms for shortest path problems in a geometric setting. Take for instance two piecewise linear paths in the Euclidean plane. The length of each path is the sum of square roots of integers (assuming the endpoint coordinates are integers). Deciding which one of the two paths is shorter is a difficult question unless the number of square roots is very small. On the other hand, deciding which one of two paths is shorter is not exactly the kind of problems that SoS was invented for.

Another problem is that algorithms employing SoS produce results for the perturbed set of objects rather than for the original ones. In certain settings, such as in computer graphics, this fact can often be ignored. However, when "unperturbed" results are needed, some postprocessing has to be performed. This paper does not deal with this issue and further work has to be done. Nevertheless, in most of the applications mentioned in this paper the postprocessing step is more or less trivial:

- In the point-in-polygon problem one can simply add a test whether or not the query point lies on a boundary edge.
- In the case of Voronoi diagrams or arrangements of hyperplanes, we identify and eliminate zero-length edges or higher-dimensional faces of zero measure.
- In the convex-hull setting, it is possible to undo the perturbation simply by merging adjacent faces if necessary; for example, in two dimensions, adjacent edges that lie on a common line, and in three dimensions, adjacent triangles contained in a common plane.

It is rather difficult, however, to use SoS or any other perturbation scheme for finding all data points on the boundary of the convex hull. This is because the perturbation may decide that a point is inside the hull if it lies on a boundary edge or face. In this case the point would be prematurely discarded. We refer to [Ya87] for a more extensive discussion of the limitations of symbolic methods aimed at resolving robustness problems in geometric algorithms.

In order to increase the credibility of our claim that SoS is indeed a practical programming tool, the second author compiled a prototype version of a SoS library [Mu88] and implemented the three-dimensional edge-skeleton algorithm of [Ed86]. We believe it is fair to say that this algorithm is an extraordinary challenge for someone who wants to do it without SoS. From run-time profiles of this program we learned that most of the computing time was spent on multiplying long integers in order to compute signs of determinants. The speed-up that we got in our implementation from replacing long-integer by normal (built-in) integer arithmetic was a factor somewhere around 10. Of course for the normal integer arithmetic to work we severely restricted the range of the coordinates that were used. In any case, this makes it clear where future work has to go if we want to produce programs that are reliable and which are as fast as software that uses floating-point arithmetic and is therefore inherently unreliable. The most promising way to eliminate this overhead factor seems to be the design of a special piece of hardware that computes the sign of determinants for integer matrices. Such effort seems justified by the versatility of determinants demonstrated in Section 5. We would like to mention, though, that even without the availability of such specialized hardware we believe that SoS is of practical value in implementing geometric algorithms. Aside from the obvious savings in time and effort for the programmer, it seems to us that the use of SoS is currently the only hope to produce geometric software that is in any sense reliable.

We end this section by pointing out a new direction for further research — it is the systematic study of primitive operations used and needed for geometric algorithms. If one undertook the venture of building a library of primitives for geometric algorithms, besides computing signs of determinants, what other operations would have to be in the collection? Is it even clear that computing the sign of a determinant is such an indispensable operation or are there less expensive ways to determine the orientation of $d+1$ points in $d$ dimensions?

# Appendix

In this appendix, we give the relevant subdeterminants, sorted in sequence of decreasing significance, needed for computing the signs of $\det \Delta_4(\varepsilon)$, $\det \Lambda_2(\varepsilon)$, $\det \Lambda_3(\varepsilon)$, and $\det \Lambda_4(\varepsilon)$. Each sequence is given in a table that also shows the corresponding $\varepsilon$-product $\varepsilon_t$ and the size $k_t$ of the matrix $M_t^{\Delta_d}$ ($M_t^{\Lambda_d}$) associated with the $(t+1)$-st significant term in the $\varepsilon$-polynomial $\det \Delta_D(\varepsilon)$ ($\det \Lambda_D(\varepsilon)$). The third column of each table shows $v_t$, the vector that encodes the subdeterminant of depth $t$. Recall that this vector was used to produce the proper sequences of subdeterminants by successive calls of procedure *Next_v*.

| $t$ | $k_t \cdot k_t$ | $v_t$ | $\det M_t^{\Lambda_2}$ | $\varepsilon_t$ |
|---|---|---|---|---|
| 0 | $2 \cdot 2$ | $[2,2;2]$ | $+\det \begin{pmatrix} \pi_{i,1} & 1 \\ \pi_{j,1} & 1 \end{pmatrix}$ | $\varepsilon()$ |
| 1 | $1 \cdot 1$ | $[1,2;2]$ | $+\det(1) = +1$ | $\varepsilon(i,1)$ |

Table -i: The 2 relevant terms of $\det \Lambda_2(\varepsilon)$.

| $t$ | $k_t \cdot k_t$ | $v_t$ | $\det M_t^{\Lambda_3}$ | $\varepsilon_t$ |
|---|---|---|---|---|
| 0 | $3 \cdot 3$ | $[3,3,3;3]$ | $+\det \begin{pmatrix} \pi_{i,1} & \pi_{i,2} & 1 \\ \pi_{j,1} & \pi_{j,2} & 1 \\ \pi_{k,1} & \pi_{k,2} & 1 \end{pmatrix}$ | $\varepsilon()$ |
| 1 | $2 \cdot 2$ | $[2,3,3;3]$ | $-\det \begin{pmatrix} \pi_{j,1} & 1 \\ \pi_{k,1} & 1 \end{pmatrix}$ | $\varepsilon(i,2)$ |
| 2 | $2 \cdot 2$ | $[1,3,3;3]$ | $+\det \begin{pmatrix} \pi_{j,2} & 1 \\ \pi_{k,2} & 1 \end{pmatrix}$ | $\varepsilon(i,1)$ |
| 3 | $2 \cdot 2$ | $[2,2,3;3]$ | $+\det \begin{pmatrix} \pi_{i,1} & 1 \\ \pi_{k,1} & 1 \end{pmatrix}$ | $\varepsilon(j,2)$ |
| 4 | $1 \cdot 1$ | $[1,2,3;3]$ | $+\det(1) = +1$ | $\varepsilon((j,2),(i,1))$ |

Table -ii: The 5 relevant terms of $\det \Lambda_3(\varepsilon)$.

| $t$ | $k_t \cdot k_t$ | $v_t$ | det $M_t^{\Lambda_4}$ | $\varepsilon_t$ |
|---|---|---|---|---|
| 0 | $4 \cdot 4$ | $[4,4,4,4;4]$ | $+\det \begin{pmatrix} \pi_{i,1} & \pi_{i,2} & \pi_{i,3} & 1 \\ \pi_{j,1} & \pi_{j,2} & \pi_{j,3} & 1 \\ \pi_{k,1} & \pi_{k,2} & \pi_{k,3} & 1 \\ \pi_{l,1} & \pi_{l,2} & \pi_{l,3} & 1 \end{pmatrix}$ | $\varepsilon()$ |
| 1 | $3 \cdot 3$ | $[3,4,4,4;4]$ | $+\det \begin{pmatrix} \pi_{j,1} & \pi_{j,2} & 1 \\ \pi_{k,1} & \pi_{k,2} & 1 \\ \pi_{l,1} & \pi_{l,2} & 1 \end{pmatrix}$ | $\varepsilon(i,3)$ |
| 2 | $3 \cdot 3$ | $[2,4,4,4;4]$ | $-\det \begin{pmatrix} \pi_{j,1} & \pi_{j,3} & 1 \\ \pi_{k,1} & \pi_{k,3} & 1 \\ \pi_{l,1} & \pi_{l,3} & 1 \end{pmatrix}$ | $\varepsilon(i,2)$ |
| 3 | $3 \cdot 3$ | $[1,4,4,4;4]$ | $+\det \begin{pmatrix} \pi_{j,2} & \pi_{j,3} & 1 \\ \pi_{k,2} & \pi_{k,3} & 1 \\ \pi_{l,2} & \pi_{l,3} & 1 \end{pmatrix}$ | $\varepsilon(i,1)$ |
| 4 | $3 \cdot 3$ | $[3,3,4,4;4]$ | $-\det \begin{pmatrix} \pi_{i,1} & \pi_{i,2} & 1 \\ \pi_{k,1} & \pi_{k,2} & 1 \\ \pi_{l,1} & \pi_{l,2} & 1 \end{pmatrix}$ | $\varepsilon(j,3)$ |
| 5 | $2 \cdot 2$ | $[2,3,4,4;4]$ | $+\det \begin{pmatrix} \pi_{k,1} & 1 \\ \pi_{l,1} & 1 \end{pmatrix}$ | $\varepsilon((j,3),(i,2))$ |
| 6 | $2 \cdot 2$ | $[1,3,4,4;4]$ | $-\det \begin{pmatrix} \pi_{k,2} & 1 \\ \pi_{l,2} & 1 \end{pmatrix}$ | $\varepsilon((j,3),(i,1))$ |
| 7 | $3 \cdot 3$ | $[2,2,4,4;4]$ | $+\det \begin{pmatrix} \pi_{i,1} & \pi_{i,3} & 1 \\ \pi_{k,1} & \pi_{k,3} & 1 \\ \pi_{l,1} & \pi_{l,3} & 1 \end{pmatrix}$ | $\varepsilon(j,2)$ |
| 8 | $2 \cdot 2$ | $[1,2,4,4;4]$ | $+\det \begin{pmatrix} \pi_{k,3} & 1 \\ \pi_{l,3} & 1 \end{pmatrix}$ | $\varepsilon((j,2),(i,1))$ |
| 9 | $3 \cdot 3$ | $[1,1,4,4;4]$ | $-\det \begin{pmatrix} \pi_{i,2} & \pi_{i,3} & 1 \\ \pi_{k,2} & \pi_{k,3} & 1 \\ \pi_{l,2} & \pi_{l,3} & 1 \end{pmatrix}$ | $\varepsilon(j,1)$ |
| 10 | $3 \cdot 3$ | $[3,3,3,4;4]$ | $+\det \begin{pmatrix} \pi_{i,1} & \pi_{i,2} & 1 \\ \pi_{j,1} & \pi_{j,2} & 1 \\ \pi_{l,1} & \pi_{l,2} & 1 \end{pmatrix}$ | $\varepsilon(k,3)$ |
| 11 | $2 \cdot 2$ | $[2,3,3,4;4]$ | $-\det \begin{pmatrix} \pi_{j,1} & 1 \\ \pi_{l,1} & 1 \end{pmatrix}$ | $\varepsilon((k,3),(i,2))$ |
| 12 | $2 \cdot 2$ | $[1,3,3,4;4]$ | $+\det \begin{pmatrix} \pi_{j,2} & 1 \\ \pi_{l,2} & 1 \end{pmatrix}$ | $\varepsilon((k,3),(i,1))$ |
| 13 | $2 \cdot 2$ | $[2,2,3,4;4]$ | $+\det \begin{pmatrix} \pi_{i,1} & 1 \\ \pi_{l,1} & 1 \end{pmatrix}$ | $\varepsilon((k,3),(j,2))$ |
| 14 | $1 \cdot 1$ | $[1,2,3,4;4]$ | $+\det(1) = +1$ | $\varepsilon((k,3),(j,2),(i,1))$ |

Table -iii: The 15 relevant terms of det $\Lambda_4(\varepsilon)$.

| $t$ | $k_t \cdot k_t$ | $v_t$ | $\det M_t^{\Delta_4}$ | $\varepsilon_t$ |
|---|---|---|---|---|
| 0 | $4 \cdot 4$ | $[5,5,5,5;5]$ | $+\det \begin{pmatrix} \pi_{i,1} & \pi_{i,2} & \pi_{i,3} & \pi_{i,4} \\ \pi_{j,1} & \pi_{j,2} & \pi_{j,3} & \pi_{j,4} \\ \pi_{k,1} & \pi_{k,2} & \pi_{k,3} & \pi_{k,4} \\ \pi_{l,1} & \pi_{l,2} & \pi_{l,3} & \pi_{l,4} \end{pmatrix}$ | $\varepsilon()$ |
| 1 | $3 \cdot 3$ | $[4,5,5,5;5]$ | $-\det \begin{pmatrix} \pi_{j,1} & \pi_{j,2} & \pi_{j,3} \\ \pi_{k,1} & \pi_{k,2} & \pi_{k,3} \\ \pi_{l,1} & \pi_{l,2} & \pi_{l,3} \end{pmatrix}$ | $\varepsilon(i,4)$ |
| 2 | $3 \cdot 3$ | $[3,5,5,5;5]$ | $+\det \begin{pmatrix} \pi_{j,1} & \pi_{j,2} & \pi_{j,4} \\ \pi_{k,1} & \pi_{k,2} & \pi_{k,4} \\ \pi_{l,1} & \pi_{l,2} & \pi_{l,4} \end{pmatrix}$ | $\varepsilon(i,3)$ |
| 3 | $3 \cdot 3$ | $[2,5,5,5;5]$ | $-\det \begin{pmatrix} \pi_{j,1} & \pi_{j,3} & \pi_{j,4} \\ \pi_{k,1} & \pi_{k,3} & \pi_{k,4} \\ \pi_{l,1} & \pi_{l,3} & \pi_{l,4} \end{pmatrix}$ | $\varepsilon(i,2)$ |
| 4 | $3 \cdot 3$ | $[1,5,5,5;5]$ | $+\det \begin{pmatrix} \pi_{j,2} & \pi_{j,3} & \pi_{j,4} \\ \pi_{k,2} & \pi_{k,3} & \pi_{k,4} \\ \pi_{l,2} & \pi_{l,3} & \pi_{l,4} \end{pmatrix}$ | $\varepsilon(i,1)$ |
| 5 | $3 \cdot 3$ | $[4,4,5,5;5]$ | $+\det \begin{pmatrix} \pi_{i,1} & \pi_{i,2} & \pi_{i,3} \\ \pi_{k,1} & \pi_{k,2} & \pi_{k,3} \\ \pi_{l,1} & \pi_{l,2} & \pi_{l,3} \end{pmatrix}$ | $\varepsilon(j,4)$ |
| 6 | $2 \cdot 2$ | $[3,4,5,5;5]$ | $+\det \begin{pmatrix} \pi_{k,1} & \pi_{k,2} \\ \pi_{l,1} & \pi_{l,2} \end{pmatrix}$ | $\varepsilon((j,4),(i,3))$ |
| 7 | $2 \cdot 2$ | $[2,4,5,5;5]$ | $-\det \begin{pmatrix} \pi_{k,1} & \pi_{k,3} \\ \pi_{l,1} & \pi_{l,3} \end{pmatrix}$ | $\varepsilon((j,4),(i,2))$ |
| 8 | $2 \cdot 2$ | $[1,4,5,5;5]$ | $+\det \begin{pmatrix} \pi_{k,2} & \pi_{k,3} \\ \pi_{l,2} & \pi_{l,3} \end{pmatrix}$ | $\varepsilon((j,4),(i,1))$ |
| 9 | $3 \cdot 3$ | $[3,3,5,5;5]$ | $-\det \begin{pmatrix} \pi_{i,1} & \pi_{i,2} & \pi_{i,4} \\ \pi_{k,1} & \pi_{k,2} & \pi_{k,4} \\ \pi_{l,1} & \pi_{l,2} & \pi_{l,4} \end{pmatrix}$ | $\varepsilon(j,3)$ |
| 10 | $2 \cdot 2$ | $[2,3,5,5;5]$ | $+\det \begin{pmatrix} \pi_{k,1} & \pi_{k,4} \\ \pi_{l,1} & \pi_{l,4} \end{pmatrix}$ | $\varepsilon((j,3),(i,2))$ |
| 11 | $2 \cdot 2$ | $[1,3,5,5;5]$ | $-\det \begin{pmatrix} \pi_{k,2} & \pi_{k,4} \\ \pi_{l,2} & \pi_{l,4} \end{pmatrix}$ | $\varepsilon((j,3),(i,1))$ |

Table -iv: The 50 relevant terms of $\det \Delta_4(\varepsilon)$ (to be continued).

| $t$ | $k_t \cdot k_t$ | $v_t$ | $\det M_t^{\Delta_4}$ | $\varepsilon_t$ |
|---|---|---|---|---|
| 12 | $3 \cdot 3$ | $[2,2,5,5;5]$ | $+\det\begin{pmatrix}\pi_{i,1} & \pi_{i,3} & \pi_{i,4}\\ \pi_{k,1} & \pi_{k,3} & \pi_{k,4}\\ \pi_{l,1} & \pi_{l,3} & \pi_{l,4}\end{pmatrix}$ | $\varepsilon(j,2)$ |
| 13 | $2 \cdot 2$ | $[1,2,5,5;5]$ | $+\det\begin{pmatrix}\pi_{k,3} & \pi_{k,4}\\ \pi_{l,3} & \pi_{l,4}\end{pmatrix}$ | $\varepsilon((j,2),(i,1))$ |
| 14 | $3 \cdot 3$ | $[1,1,5,5;5]$ | $-\det\begin{pmatrix}\pi_{i,2} & \pi_{i,3} & \pi_{i,4}\\ \pi_{k,2} & \pi_{k,3} & \pi_{k,4}\\ \pi_{l,2} & \pi_{l,3} & \pi_{l,4}\end{pmatrix}$ | $\varepsilon(j,1)$ |
| 15 | $3 \cdot 3$ | $[4,4,4,5;5]$ | $-\det\begin{pmatrix}\pi_{i,1} & \pi_{i,2} & \pi_{i,3}\\ \pi_{j,1} & \pi_{j,2} & \pi_{j,3}\\ \pi_{l,1} & \pi_{l,2} & \pi_{l,3}\end{pmatrix}$ | $\varepsilon(k,4)$ |
| 16 | $2 \cdot 2$ | $[3,4,4,5;5]$ | $-\det\begin{pmatrix}\pi_{j,1} & \pi_{j,2}\\ \pi_{l,1} & \pi_{l,2}\end{pmatrix}$ | $\varepsilon((k,4),(i,3))$ |
| 17 | $2 \cdot 2$ | $[2,4,4,5;5]$ | $+\det\begin{pmatrix}\pi_{j,1} & \pi_{j,3}\\ \pi_{l,1} & \pi_{l,3}\end{pmatrix}$ | $\varepsilon((k,4),(i,2))$ |
| 18 | $2 \cdot 2$ | $[1,4,4,5;5]$ | $-\det\begin{pmatrix}\pi_{j,2} & \pi_{j,3}\\ \pi_{l,2} & \pi_{l,3}\end{pmatrix}$ | $\varepsilon((k,4),(i,1))$ |
| 19 | $2 \cdot 2$ | $[3,3,4,5;5]$ | $+\det\begin{pmatrix}\pi_{i,1} & \pi_{i,2}\\ \pi_{l,1} & \pi_{l,2}\end{pmatrix}$ | $\varepsilon((k,4),(j,3))$ |
| 20 | $1 \cdot 1$ | $[2,3,4,5;5]$ | $-\det(\pi_{l,1}) = -\pi_{l,1}$ | $\varepsilon((k,4),(j,3),(i,2))$ |
| 21 | $1 \cdot 1$ | $[1,3,4,5;5]$ | $+\det(\pi_{l,2}) = +\pi_{l,2}$ | $\varepsilon((k,4),(j,3),(i,1))$ |
| 22 | $2 \cdot 2$ | $[2,2,4,5;5]$ | $-\det\begin{pmatrix}\pi_{i,1} & \pi_{i,3}\\ \pi_{l,1} & \pi_{l,3}\end{pmatrix}$ | $\varepsilon((k,4),(j,2))$ |
| 23 | $1 \cdot 1$ | $[1,2,4,5;5]$ | $-\det(\pi_{l,3}) = -\pi_{l,3}$ | $\varepsilon((k,4),(j,2),(i,1))$ |
| 24 | $2 \cdot 2$ | $[1,1,4,5;5]$ | $+\det\begin{pmatrix}\pi_{i,2} & \pi_{i,3}\\ \pi_{l,2} & \pi_{l,3}\end{pmatrix}$ | $\varepsilon((k,4),(j,1))$ |
| 25 | $3 \cdot 3$ | $[3,3,3,5;5]$ | $+\det\begin{pmatrix}\pi_{i,1} & \pi_{i,2} & \pi_{i,4}\\ \pi_{j,1} & \pi_{j,2} & \pi_{j,4}\\ \pi_{l,1} & \pi_{l,2} & \pi_{l,4}\end{pmatrix}$ | $\varepsilon(k,3)$ |
| 26 | $2 \cdot 2$ | $[2,3,3,5;5]$ | $-\det\begin{pmatrix}\pi_{j,1} & \pi_{j,4}\\ \pi_{l,1} & \pi_{l,4}\end{pmatrix}$ | $\varepsilon((k,3),(i,2))$ |
| 27 | $2 \cdot 2$ | $[1,3,3,5;5]$ | $+\det\begin{pmatrix}\pi_{j,2} & \pi_{j,4}\\ \pi_{l,2} & \pi_{l,4}\end{pmatrix}$ | $\varepsilon((k,3),(i,1))$ |
| 28 | $2 \cdot 2$ | $[2,2,3,5;5]$ | $+\det\begin{pmatrix}\pi_{i,1} & \pi_{i,4}\\ \pi_{l,1} & \pi_{l,4}\end{pmatrix}$ | $\varepsilon((k,3),(j,2))$ |
| 29 | $1 \cdot 1$ | $[1,2,3,5;5]$ | $+\det(\pi_{l,4}) = +\pi_{l,4}$ | $\varepsilon((k,3),(j,2),(i,1))$ |

Table 6 continued.

| $t$ | $k_t \cdot k_t$ | $v_t$ | $\det M_t^{\Delta_4}$ | $\varepsilon_t$ |
|---|---|---|---|---|
| 30 | $2 \cdot 2$ | $[1,1,3,5;5]$ | $-\det \begin{pmatrix} \pi_{i,2} & \pi_{i,4} \\ \pi_{l,2} & \pi_{l,4} \end{pmatrix}$ | $\varepsilon((k,3),(j,1))$ |
| 31 | $3 \cdot 3$ | $[2,2,2,5;5]$ | $-\det \begin{pmatrix} \pi_{i,1} & \pi_{i,3} & \pi_{i,4} \\ \pi_{j,1} & \pi_{j,3} & \pi_{j,4} \\ \pi_{l,1} & \pi_{l,3} & \pi_{l,4} \end{pmatrix}$ | $\varepsilon(k,2)$ |
| 32 | $2 \cdot 2$ | $[1,2,2,5;5]$ | $-\det \begin{pmatrix} \pi_{j,3} & \pi_{j,4} \\ \pi_{l,3} & \pi_{l,4} \end{pmatrix}$ | $\varepsilon((k,2),(i,1))$ |
| 33 | $2 \cdot 2$ | $[1,1,2,5;5]$ | $+\det \begin{pmatrix} \pi_{i,3} & \pi_{i,4} \\ \pi_{l,3} & \pi_{l,4} \end{pmatrix}$ | $\varepsilon((k,2),(j,1))$ |
| 34 | $3 \cdot 3$ | $[1,1,1,5;5]$ | $+\det \begin{pmatrix} \pi_{i,2} & \pi_{i,3} & \pi_{i,4} \\ \pi_{j,2} & \pi_{j,3} & \pi_{j,4} \\ \pi_{l,2} & \pi_{l,3} & \pi_{l,4} \end{pmatrix}$ | $\varepsilon(k,1)$ |
| 35 | $3 \cdot 3$ | $[4,4,4,4;5]$ | $+\det \begin{pmatrix} \pi_{i,1} & \pi_{i,2} & \pi_{i,3} \\ \pi_{j,1} & \pi_{j,2} & \pi_{j,3} \\ \pi_{k,1} & \pi_{k,2} & \pi_{k,3} \end{pmatrix}$ | $\varepsilon(l,4)$ |
| 36 | $2 \cdot 2$ | $[3,4,4,4;5]$ | $+\det \begin{pmatrix} \pi_{j,1} & \pi_{j,2} \\ \pi_{k,1} & \pi_{k,2} \end{pmatrix}$ | $\varepsilon((l,4),(i,3))$ |
| 37 | $2 \cdot 2$ | $[2,4,4,4;5]$ | $-\det \begin{pmatrix} \pi_{j,1} & \pi_{j,3} \\ \pi_{k,1} & \pi_{k,3} \end{pmatrix}$ | $\varepsilon((l,4),(i,2))$ |
| 38 | $2 \cdot 2$ | $[1,4,4,4;5]$ | $+\det \begin{pmatrix} \pi_{j,2} & \pi_{j,3} \\ \pi_{k,2} & \pi_{k,3} \end{pmatrix}$ | $\varepsilon((l,4),(i,1))$ |
| 39 | $2 \cdot 2$ | $[3,3,4,4;5]$ | $-\det \begin{pmatrix} \pi_{i,1} & \pi_{i,2} \\ \pi_{k,1} & \pi_{k,2} \end{pmatrix}$ | $\varepsilon((l,4),(j,3))$ |
| 40 | $1 \cdot 1$ | $[2,3,4,4;5]$ | $+\det(\pi_{k,1}) = +\pi_{k,1}$ | $\varepsilon((l,4),(j,3),(i,2))$ |
| 41 | $1 \cdot 1$ | $[1,3,4,4;5]$ | $-\det(\pi_{k,2}) = -\pi_{k,2}$ | $\varepsilon((l,4),(j,3),(i,1))$ |
| 42 | $2 \cdot 2$ | $[2,2,4,4;5]$ | $+\det \begin{pmatrix} \pi_{i,1} & \pi_{i,3} \\ \pi_{k,1} & \pi_{k,3} \end{pmatrix}$ | $\varepsilon((l,4),(j,2))$ |
| 43 | $1 \cdot 1$ | $[1,2,4,4;5]$ | $+\det(\pi_{k,3}) = +\pi_{k,3}$ | $\varepsilon((l,4),(j,2),(i,1))$ |
| 44 | $2 \cdot 2$ | $[1,1,4,4;5]$ | $-\det \begin{pmatrix} \pi_{i,2} & \pi_{i,3} \\ \pi_{k,2} & \pi_{k,3} \end{pmatrix}$ | $\varepsilon((l,4),(j,1))$ |
| 45 | $2 \cdot 2$ | $[3,3,3,4;5]$ | $+\det \begin{pmatrix} \pi_{i,1} & \pi_{i,2} \\ \pi_{j,1} & \pi_{j,2} \end{pmatrix}$ | $\varepsilon((l,4),(k,3))$ |
| 46 | $1 \cdot 1$ | $[2,3,3,4;5]$ | $-\det(\pi_{j,1}) = -\pi_{j,1}$ | $\varepsilon((l,4),(k,3),(i,2))$ |
| 47 | $1 \cdot 1$ | $[1,3,3,4;5]$ | $+\det(\pi_{j,2}) = +\pi_{j,2}$ | $\varepsilon((l,4),(k,3),(i,1))$ |
| 48 | $1 \cdot 1$ | $[2,2,3,4;5]$ | $+\det(\pi_{i,1}) = +\pi_{i,1}$ | $\varepsilon((l,4),(k,3),(j,2))$ |
| 49 | $0 \cdot 0$ | $[1,2,3,4;5]$ | $+\det() = +1$ | $\varepsilon((l,4),(k,3),(j,2),(i,1))$ |

Table 6 continued.